\documentclass[ijoc,nonblindrev,12pt]{informs3} 

\OneAndAHalfSpacedXII 

\usepackage{ booktabs, makecell, multirow, multirow, float,subfigure}
\usepackage[colorlinks, linkcolor=red, anchorcolor=blue, citecolor=blue]{hyperref}
\usepackage{import}
\usepackage{bbm}
\usepackage{color, xcolor}

\usepackage{graphicx}
\usepackage{amsmath}
\usepackage{amsfonts}

\usepackage{bm,todonotes}
\usepackage{comment}

\newcommand{\cO}{\mathcal{O}}

\newcommand{\EE}{\mathbb{E}}
\newcommand{\RR}{\mathbb{R}}

\newcommand{\xb}{\mathbf{x}}

\newcommand{\ub}{\mathbf{u}}

\newcommand{\I}{\mathbf{I}}

\usepackage{bbm}

\newcommand{\ti}[1]{\tilde{#1}}
\def\diag{\mathrm{diag}}
\newcommand{\norm}[1]{\left\|#1\right\|}
\newcommand{\dotprod}[1]{\left\langle #1\right\rangle}
\def\tr{\mathrm{tr}}

\usepackage{amsmath}

\usepackage{multirow}
\usepackage{natbib}
\usepackage{algorithm}
\usepackage{algorithmic}

 \bibpunct[, ]{(}{)}{,}{a}{}{,}%
 \def\bibfont{\small}%
\newcommand{\mL}{\widehat{L}}

\usepackage{bbm}

\newcommand{\x}{\mathbf{x}}
\newcommand{\y}{\mathbf{y}}
\newcommand{\g}{\mathbf{g}}
\newcommand{\s}{\mathbf{s}}
\newcommand{\bv}{\mathbf{v}}
\newcommand{\bw}{\mathbf{w}}

\usepackage{comment}
\usepackage{epstopdf}
\usepackage[colorlinks,
            linkcolor=red,
            anchorcolor=blue,
            citecolor=blue
            ]{hyperref}
\usepackage{natbib}

\newcommand{\xiangyu}[1]{{\textbf{\textcolor{red}{[Xiangyu: #1]}}}}

\def\diag{\mathrm{diag}}

\def\tr{\mathrm{tr}}
\newcommand{\sti}{\mathbf{s}_t^{(i)}}
\newcommand{\si}{\mathbf{s}^{(i)}}
\newcommand{\tg}{\tilde{\g}}
\newcommand{\tv}{\tilde{\bv}}

\TheoremsNumberedThrough     

\EquationsNumberedThrough    

\newcommand{\haishan}[1]{{\textbf{\textcolor{orange}{[haishan: #1]}}}}
\begin{document}

\TITLE{A Unified Zeroth-Order Optimization Framework via Oblivious Randomized Sketching}

\ARTICLEAUTHORS{

\AUTHOR{Haishan Ye}
\AFF{
	School of Management, Xi'an Jiaotong University, Xi'an 710049, China  \\ \EMAIL{yehaishan@xjtu.edu.cn}\URL{}}

\AUTHOR{Xiangyu Chang}
\AFF{
		School of Management, Xi'an Jiaotong University, Xi'an 710049, China \\ \EMAIL{xiangyuchang@xjtu.edu.cn}}

\AUTHOR{Xi Chen}
\AFF{
		Stern School of Business, New York University, New York, NY 10012, USA\\ \EMAIL{xc13@stern.nyu.edu}}

}

\ABSTRACT{%
We propose a new framework for analyzing zeroth-order optimization (ZOO) from the perspective of \emph{oblivious randomized sketching}.
In this framework, commonly used gradient estimators in ZOO—such as finite difference (FD) and random finite difference (RFD)—are unified through a general sketch-based formulation.
By introducing the concept of oblivious randomized sketching, we show that properly chosen sketch matrices can significantly reduce the high variance of RFD estimates and enable \emph{high-probability} convergence guarantees of ZOO, which are rarely available in existing RFD analyses.

\noindent We instantiate the framework on convex quadratic objectives and derive a query complexity of $\tilde{\mathcal{O}}(\mathrm{tr}(A)/L \cdot L/\mu\log\frac{1}{\epsilon})$ to achieve a $\epsilon$-suboptimal solution, where $A$ is the Hessian, $L$ is the largest eigenvalue of $A$, and $\mu$ denotes the strong convexity parameter.
This complexity can be substantially smaller than the standard query complexity of ${\cO}(d\cdot L/\mu \log\frac{1}{\epsilon})$ that is linearly dependent on problem dimensionality, especially when $A$ has rapidly decaying eigenvalues.
These advantages naturally extend to more general settings, including strongly convex and Hessian-aware optimization.

\noindent  Overall, this work offers a novel sketch-based perspective on ZOO that explains why and when RFD-type methods can achieve \emph{weakly dimension-independent} convergence in general smooth problems, providing both theoretical foundations and practical implications for ZOO.
}


\KEYWORDS{Zeroth-order Optimization, Gradient Estimation, Randomized Sketching, Trace Estimation}
\HISTORY{\today}

\maketitle

%

\section{Introduction}\label{sec:Intro}

This paper considers a specific unconstrained optimization problem of the form
\begin{equation}\label{Eq: main objective}
	\min_{\x\in\RR^n} \phi(\x),
\end{equation}
where we assume the noisy approximation $f(\x) = \phi(\x) + \zeta(\x)$ is accessible ($\zeta(\x)$ describes the noisy term), while $\phi(\x)$ or $\nabla \phi(\x)$ are not available.
In other words, we aim to minimize $\phi(\x)$ using only the information of its noisy approximation $f(\x)$.
The incidence of the above problem \eqref{Eq: main objective} has shown marked escalation in modern computational applications, including but not limited to 
memory-efficient fine-tuning of large language models \citep{malladi2023fine}, black-box adversarial attacks \citep{ilyas2018black}, reinforcement learning \citep{choromanski2018structured}, optimization with bandit feedback \citep{bubeck2012regret,wang2025adaptivity}, and distributionally robust optimization \citep{lam2024distributionally}.

To address problem \eqref{Eq: main objective},
the so-called zeroth-order optimization (ZOO) methods have been extensively investigated, which operate without explicit access to gradient information and rely solely on function evaluations. 
The fundamental idea of ZOO is to compute an estimate of the gradient of $\nabla\phi(\x)$ at the point $\x$ by $f(\x)$, and then adopt a gradient-based method to solve \eqref{Eq: main objective}~\citep{ghadimi2013stochastic,nesterov2017random}.
Therefore, the most critical component of ZOO algorithms is the construction of accurate and efficient gradient estimates from only function value queries of $f(\x)$.

Among the various gradient estimation techniques in ZOO, two classical approaches are widely used. 
The first is the \textit{finite difference} (FD) method, which approximates the gradient by evaluating the function at points sampled along coordinate directions. The most common approach is the classical Kiefer-Wolfowitz scheme~\citep{kiefer1952stochastic}, where the function is evaluated at a single nearby point along each of the $d$ coordinate axes to estimate each partial derivative.
A more accurate variant, known as the \textit{central finite difference}, samples function values at two points symmetrically located around the current iterate along each coordinate, thereby improving the approximation accuracy at the cost of additional function evaluations.
The second popular approach is the \textit{random finite difference} (RFD) method, which constructs gradient estimates by aggregating directional derivative approximations along randomly chosen directions. 
Prominent examples of RFD techniques include the Simultaneous Perturbation Stochastic Approximation (SPSA) method introduced by~\citet{spall1992multivariate}, which employs perturbation distributions characterized by significant mass at zero. Another well-studied class is Gaussian smoothing approaches~\citep{ghadimi2013stochastic,nesterov2017random}, which utilize multivariate Gaussian distributions to generate perturbation vectors. 
More recent research has explored the generation of dependent perturbation vectors, such as those constructed via orthogonal sampling schemes~\citep{choromanski2018structured} and Dirichlet distribution~\citep{lam2024distributionally}. 
Comprehensive empirical and theoretical comparisons of these gradient estimation methods are available in \citet{berahas2022theoretical}.

FD and RFD methods are widely employed in the ZOO. 
On one hand, FD directly approximates gradients through deterministic perturbations along each coordinate direction. 
Consequently, FD inherits many desirable properties from first-order methods, including convergence behaviors and compatibility with established acceleration techniques such as quasi-Newton methods, significantly enhancing practical optimization performance~\citep{ye2021approximate}. 
However, the computational burden of FD scales linearly with dimensionality $d$, as it requires evaluating all $d$ partial derivatives at each iteration. 
Thus, FD methods quickly become computationally infeasible for high-dimensional ZOO problems. 
Random coordinate descent, which can be regarded as a kind of randomized variant of FD, can only require $\cO(1)$ samples compared to FD. 
However, random coordinate descent can also only achieve a query complexity linearly dependent on the dimension $d$~\citep{nesterov2012efficiency}.
\citet{nesterov2012efficiency} propose to apply importance sampling to random coordinate descent and have shown theoretical promise in overcoming the dimensionality bottleneck, achieving convergence rates independent of $d$ under certain conditions.
Nevertheless, the practical deployment of importance sampling remains challenging. 
Specifically, optimal importance sampling distributions must be carefully constructed based on smoothness parameters for each coordinate, necessitating computationally expensive parameter estimation steps \citep{nesterov2012efficiency}. 
This dependence on dimensionally scaled parameter estimations thus reintroduces dimensional constraints into practical implementations.

On the other hand, RFD methods estimate gradients by randomly sampling directions, drastically reducing per-iteration computational complexity compared to FD in practice. 
However, current theory shows that RFD still achieves a query complexity linearly dependent on dimension~\citep{ghadimi2013stochastic,nesterov2017random}.
In deed, RFD continues to be the predominant choice in high-dimensional black-box optimization applications, such as fine-tuning large language models~\citep{malladi2023fine} and adversarial attacks~\citep{ilyas2018black}. 
\textit{This prevalent adoption highlights a notable discrepancy between theoretical limitations and practical usage: while RFD methods are empirically favored due to lower computational overhead per iteration, their theoretical efficacy, particularly concerning convergence and accuracy, is often not significantly superior to FD.}

Acknowledging this fundamental gap, \citet{scheinberg2022finite} explicitly formulated several open questions to clarify the efficacy and limitations of RFD-based methods in high-dimensional settings. 
Among these, two critical questions stand out prominently:
\begin{itemize}
    \item Given extremely high-dimensional problems (e.g., $d \approx 10^6$) where even a single full gradient evaluation using FD is computationally prohibitive, can RFD-based gradient descent algorithms reliably achieve meaningful progress within fewer than $\mathcal{O}(d)$ iterations, despite inherently high gradient estimation variance?
    \item Line search methods based on FD can leverage quasi-Newton techniques to approximate Hessians effectively. 
    Conversely, implementing similar quasi-Newton methods with RFD estimates is problematic due to their inherently low accuracy. 
    Is it possible, therefore, to reconcile the accuracy gap between RFD and FD approaches in the context of second-order optimization schemes?
\end{itemize}

To systematically address the above open questions, we propose a unified framework that integrates oblivious randomized sketching techniques into the ZOO. 
Our approach is motivated by the observation that gradient estimations using either FD or RFD can be unified under the representation involving a pre-defined matrix $S$ (see Lemma \ref{lem:g_prop}). 
By selecting $S$ explicitly as a random sketching matrix, we can effectively control the high variance typically encountered in RFD gradient estimation.
Additionally, this randomized sketching approach significantly alleviates the theoretical dimension-dependency in the convergence rate. 
Specifically, our framework achieves a high-probability convergence rate for $\epsilon$ precision with query complexity $\Tilde{\cO}\left(\frac{\tr(A)}{\mu}\log\frac{1}{\epsilon}\right)$, where $A$ is related to the second-order information of the objective (see Assumption \ref{ass:phi}), $\mu$ is the strong convexity parameter (see Assumption \ref{ass:mu}) and $\Tilde{\cO}(\cdot)$ hides constant and $\log\log\frac{1}{\epsilon}$ terms.
This quantity  $\Tilde{\cO}\left(\frac{\tr(A)}{\mu}\log\frac{1}{\epsilon}\right)$ can be substantially smaller than the traditional $\cO(d\cdot L/\mu \cdot \log\frac{1}{\epsilon})$ complexity encountered in many practical scenarios~(See the Remark \ref{remark_1} and empirical analysis in Section \ref{sec:experiment}).
Moreover, by extending this unified framework through approximate Hessian estimation~\citep{hansen2016cma,ye2018hessian}, we propose a Hessian-aware zeroth-order optimization framework with randomized sketching. This extended approach not only retains the advantages of randomized sketching but also exhibits enhanced convergence properties. 



\subsection{Related Work}

From the classical theory of convex optimization, an accurate gradient estimation is crucial for achieving a fast convergence rate. 
Consequently, a large body of work in ZOO focuses on constructing accurate gradient estimates via zeroth-order oracles, typically under the Frobenius norm metric, using FD or RFD methods. 
The goal is to find an approximate gradient $\g(\x)$ that minimizes the error $\norm{\g(\x) - \nabla \phi(\x)}$.  

A representative FD method is the classical Kiefer--Wolfowitz scheme, which estimates $d$ partial derivatives along each coordinate direction using zeroth-order oracles \citep{kiefer1952stochastic}. 
This approach yields a high-precision gradient approximation but requires $2d$ oracle calls per iteration. 
Combined with the convergence theory of gradient descent \citep{nesterov2003introductory}, the query complexity of the Kiefer--Wolfowitz method scales linearly with the dimension $d$.  
An alternative is the linear interpolation method \citep{powell1972unconstrained,conn1997convergence,wild2008orbit,conn2009introduction}, which also requires $\cO(d)$ samples in a neighborhood of $\x$. 
Unlike the Kiefer--Wolfowitz method, interpolation employs $d$ linearly independent search directions instead of coordinate directions. 
Although this method can also achieve high-precision gradient approximations, its query complexity still exhibits linear dependence on the problem dimension due to the $\cO(d)$ oracle calls required per iteration.

Recently, RFD methods with different random perturbations have gained popularity in computational and scientific applications \citep{ilyas2018black,malladi2023fine,wang2025adaptivity}. 
Typical perturbation distributions include the Gaussian \citep{malladi2023fine,nesterov2017random,ghadimi2013stochastic} and the spherical distribution \citep{zhang2024dpzero}. 
However, when aiming for a gradient estimation with small error $\norm{\g(\x) - \nabla \phi(\x)}$, \citet{berahas2022theoretical} show that RFD with either Gaussian or spherical perturbations requires $\cO(d)$ samples per iteration, resulting in query complexity that scales linearly with the problem dimension.  

On the other hand, practical implementations of RFD typically use only $\cO(1)$ samples per iteration. 
In this regime, the variance of the estimated gradient $\g(\x)$, namely $\EE\norm{\g(\x)}^2$, dominates the convergence behavior of zeroth-order algorithms \citep{nesterov2017random,ghadimi2013stochastic}. 
\citet{nesterov2017random} demonstrate that for RFD with Gaussian search directions, the variance scales as $\cO(d\cdot \norm{\nabla \phi(\x)}^2)$ when the perturbation is sufficiently small. 
Consequently, the query complexity of Gaussian RFD is also linearly dependent on $d$. 
More recently, \citet{ma2025revisiting} show that under the condition $\EE[\s\s^\top] = I_d$, where $\s$ denotes the random search direction, the RFD method using spherical directions scaled by $\sqrt{d}$ achieves the minimum variance $\EE\norm{\g(\x) - \nabla \phi(\x)}^2$. 
Nevertheless, this minimum variance remains $\cO(d\cdot \norm{\nabla \phi(\x)}^2)$, which still results in a query complexity that grows linearly with the dimension.  

The breakthrough in reducing the dimension dependence of RFD methods lies in the work of \citet{yue2023zeroth}. 
They first analyze the quadratic convex problem, in which the Hessian is a constant matrix $A$. 
In this setting, they study the variance of the gradient estimate in the $A$-norm, that is, $\EE\norm{\g(\x)}_A^2$. 
They prove that for random Gaussian perturbations, the variance satisfies $\EE\norm{\g(\x)}_A^2 = \cO\!\left(\tr(A)\cdot \norm{\nabla \phi(\x)}^2\right)$. 
This variance bound immediately implies that Gaussian RFD achieves a query complexity depending on $\frac{\tr(A)}{\lambda_{\max}(A)}$, where $\lambda_{\max}(A)$ denotes the largest eigenvalue of $A$. 
When $A$ exhibits fast eigenvalue decay, the ratio $\frac{\tr(A)}{\lambda_{\max}(A)}$ can be much smaller than $d$, thereby yielding only a weak dependence on the dimension.  
It is worth emphasizing, however, that the analysis of \citet{yue2023zeroth} fundamentally relies on properties specific to the Gaussian distribution. 
Whether analogous results can be established for other perturbation distributions without exploiting such Gaussian-specific structures remains an open question.  
Almost simultaneously, \citet{malladi2023fine} obtain nearly the same result to explain why ZOO can be effective in fine-tuning large language models. 
More recently, \citet{zhang2024dpzero} extend this line of research by showing similar results for spherical perturbations. 
For comparison, \citet{nesterov2012efficiency} achieves comparable query complexity in the context of random coordinate descent with importance sampling, but their method requires explicit access to all diagonal entries of the Hessian $A$, which is impractical in many high-dimensional applications.

Although several studies have shown that RFD with Gaussian or spherical perturbations can achieve weak dimension dependency, their analyses heavily rely on the special properties of these distributions, such as rotational invariance and exact moment conditions. 
As a result, it is still unclear what general sufficient conditions on the random search direction are required to ensure weak dimension dependency (i.e., less than $\cO(d))$). 
Moreover, most existing results are established only in expectation, and how to obtain high-probability guarantees for RFD methods remains a problem worthy of study.  
This gap directly relates to the open questions that have been discussed in the introduction section proposed by \citet{scheinberg2022finite}.
In this paper, we provide an affirmative answer. 
We provide a unified framework and identify a set of sufficient conditions on the distribution of random search directions that guarantee weak dimension dependency of RFD methods. 
Furthermore, under these conditions, we establish high-probability results, thereby systematically resolving the open question raised in \citet{scheinberg2022finite} and advancing the theory of ZOO methods beyond the Gaussian and spherical cases.  

\subsection{Contributions}
The significance of our work is threefold.

\begin{itemize}
    \item \textbf{Indicating connection between oblivious randomized sketching matrices and variance reduction in RFD.} 
    We establish a unified theoretical framework that bridges the oblivious randomized sketching techniques and RFD methods. 
    By interpreting RFD gradient estimates through the lens of sketching, we reveal that oblivious randomized sketching matrices can effectively control the variance of gradient estimates.
    This insight offers a principled approach to designing low-variance, {weakly dimension-independent} zeroth-order algorithms using random search directions that extend beyond Gaussian and spherical distributions.
    
    \item \textbf{Providing high-probability convergence guarantees with weakly dimension-independent query complexity.} 
    Leveraging the above connection, we analyze our zeroth-order algorithm in both quadratic convex and strongly convex settings. 
    Our analysis yields high-probability convergence bounds whose query complexity depends only on the intrinsic geometry of the problem, specifically $\mathrm{tr}(A)/\mu$ where $A$ is related to the Hessian information and $\mu$ is the strong convexity parameter.
    This significantly improves upon the dependence on $d \cdot \frac{L}{\mu}$ of the query complexity of classical RFD methods in many theoretical and practical situations and positively resolves the open question raised by \citet{scheinberg2022finite}.

    \item \textbf{Extending to the Hessian-aware case with adaptive step size via trace estimation.} 
    We further extend our framework to the general Hessian-aware setting, where the approximate Newton direction is used instead of the standard gradient direction.
    We show that the same sketching-based framework can be applied to efficiently estimate the trace of the Hessian, enabling automatic step size selection in practice.
    This makes our algorithm more adaptive and scalable to high-dimensional problems, while retaining favorable theoretical guarantees.
\end{itemize}

Through these theoretical advancements and practical algorithmic enhancements, our research systematically addresses and provides concrete answers to the open questions posed by \cite{scheinberg2022finite}.

\section{Problem Setting and Challenges}\label{sec:problem_setting}

\subsection{Notations}

We introduce essential notations used throughout this paper. 
A symmetric matrix $A \in \RR^{d \times d}$ is called positive definite if it holds that $\x^\top A \x > 0$ for any non-zero vector $\x \in \RR^d$.
Let $ A = U \Lambda U^\top $ be the eigenvalue decomposition of positive definite matrix $A$, where $U$ is an orthonormal matrix which contains the eigenvectors of $A$ and $\Lambda = \diag(\lambda_1,\lambda_2,\dots, \lambda_d)$ is a diagonal matrix with $\lambda_1 \ge \lambda_2\ge\dots\ge \lambda_d >0$.
Using the eigenvalues, we can define the condition number of $A$ as $\kappa(A) \triangleq \frac{\lambda_1}{\lambda_d}$. 
We can also define the square root of a positive definite matrix $A$ as $A^{1/2} = U\Lambda^{1/2}U^\top$.
This paper uses $\tr(A) = \sum_{i=1}^d A_{ii}$ to denote the trace of the matrix $A$.
Given two vectors $\x, \y\in\RR^d$, we use $\dotprod{\x,\y} = \sum_{i=1}^{d} x_i y_i$ to denote the inner product of $\x$ and $\y$.
Let $\{\mathbf{e}_i\}_{i=1}^d$ be the standard basis in $\mathbb{R}^d$ with the $i$-th coordinate as 1 and the other coordinates as 0.

Additionally, $\|A\|_{F} \triangleq (\sum_{i,j}A_{ij}^{2})^{1/2}=(\sum_{i}\lambda_{i}^{2})^{1/2}$
is the Frobenius norm of a positive definite $A$ and
$\|A\|_2\triangleq \lambda_{1}$ is the spectral norm. 
For a vector $\x\in\RR^d$, we define its norm as $\norm{\x} \triangleq \sqrt{\sum_{i=1}^d x_i^2}$.
Given a symmetric positive definite matrix $A$, $\|\x\|_A \triangleq \|A^{1/2}\x\|$ is called the $A$-norm of $\x$.
Given square matrices $A$ and $B$ with the same size, we denote $A \preceq B$ if $B-A$ is positive semidefinite.

\subsection{Assumptions}

\begin{assumption} \label{ass:phi}
	The objective function $\phi(\x)$ is convex quadratic with the Hessian matrix $A \in\RR^{d\times d}$. Accordingly, it holds that
	\begin{equation}\label{eq:phi}
		\phi(\y) = \phi(\x) + \dotprod{\nabla \phi(\x), \y - \x} + \frac{1}{2} \norm{\y - \x}_A^2,
	\end{equation}
    where $0 \preceq  A\preceq L\cdot I$ and $L>0$ is a constant.
\end{assumption}
\begin{assumption}\label{ass:mu}
	The objective function $\phi(\x)$ is $\mu$-strongly convex. Then it holds that
	\begin{equation}\label{eq:mu}
		\phi(\y) \ge \phi(\x) + \dotprod{\nabla \phi(\x), \y - \x} + \frac{\mu}{2} \norm{\y - \x}^2.
	\end{equation}
\end{assumption}

\begin{assumption}\label{ass:zeta}
	There is a constant $\sigma\ge 0$, such that $|f(\x) - \phi(\x)| \le \sigma$ for all $\x \in \RR^d$. 
\end{assumption}

\begin{assumption}\label{ass:L}
	The objective function $\phi(\x)$ is $L$-smooth. That is, for any $\x,\y\in\RR^d$, it holds 
	\begin{equation}
		\phi(\y) \leq \phi(\x) + \dotprod{\nabla \phi(\x), \y - \x} + \frac{L}{2}\norm{\y - \x}^2.
	\end{equation}
\end{assumption}
\begin{assumption}\label{ass:mL}
	The objective function $\phi(\x)$ admits an $\mL$-Lipschitz continuous Hessian. That is, for any $\x,\y\in\RR^d$, it holds 
	\begin{equation}
		\norm{\nabla^2 \phi(\x) - \nabla^2\phi(\y)}_2 \le \mL \cdot \norm{\x - \y}.
	\end{equation}
\end{assumption}
Assumption~\ref{ass:phi}-Assumption~\ref{ass:mu} and Assumption~\ref{ass:L}-Assumption~\ref{ass:mL} are about the properties of the objective function. These assumptions are widely used in the convergence analysis in convex analysis \citep{nesterov2003introductory}.
Assumption~\ref{ass:zeta} is related to the noise.
Assumption~\ref{ass:zeta} is widely used in ZOO which assumes that the noise is bounded \citep{berahas2022theoretical,yue2023lower}.

We first analyze the quadratic surrogate setting in Section \ref{sec:convergence_convex_quad} under Assumption~\ref{ass:phi}.
This assumption is standard and representative in the ZOO literature and has been explicitly adopted in recent work on weak dimension dependency for zeroth-order methods~\citep{yue2023zeroth}.
Building on this baseline, Sections~\ref{sec:hessian_aware} and~\ref{sec:convegence_smooth_convex} extend our results along two axes that require additional regularity: (i) \emph{smooth and strong convexity} (Assumption~\ref{ass:mu} and \ref{ass:L}); and
(ii) \emph{Hessian-aware Problem} (Assumptions~\ref{ass:mL}).
Finally, Assumption~\ref{ass:zeta} imposes a bounded-noise property on the function oracle, a common condition in ZOO analyses~\citep{berahas2022theoretical,scheinberg2022finite}.

\subsection{Gradient Estimation and Its Challenges}

Let $S \in \mathbb{R}^{d \times \ell}$ be a pre-defined matrix, and let $\si$ denote its $i$-th column. 
We approximate the gradient of $\phi(\mathbf{x})$ at a given point $\mathbf{x}$ by
\begin{equation}
\mathbf{g}(\mathbf{x}) = \sum_{i=1}^\ell \frac{f(\mathbf{x}+\alpha \si) - f(\mathbf{x}-\alpha \si)}{2\alpha} \si. \label{eq:g}
\end{equation}
Notably, when columns $\si$ of $S$ correspond to the $i$-th standard basis vectors of $\{\mathbf{e}_i\}_{i=1}^d$, $\ell=d$ and $S=I_d$ ($I_d$ is the identical matrix), $\mathbf{g}(\mathbf{x})$ reduces to the classical FD approximation associated with the Kiefer-Wolfowitz scheme. 
Alternatively, if each column $\si$ is randomly sampled from either the standard basis or a given probability distribution (e.g., a multivariate Gaussian distribution), then the resulting $\mathbf{g}(\mathbf{x})$ constitutes an RFD approximation of $\nabla\phi(\mathbf{x})$.

The subsequent lemma rigorously establishes a clear and explicit relationship between the gradient estimate $\mathbf{g}(\mathbf{x})$, the pre-defined matrix $S$, and the true gradient $\nabla\phi(\mathbf{x})$ under Assumption \ref{ass:phi}. 
This characterization further elucidates why selecting $S$ as a randomized sketching matrix can effectively control the variance in gradient estimates and improve algorithmic performance in the ZOO.

\begin{lemma}\label{lem:g_prop}
	Let the objective function $\phi(\x)$ satisfy Assumption~\ref{ass:phi} and Assumption~\ref{ass:zeta}. 
	We assume that  $S\in\RR^{d \times \ell}$ is a pre-defined matrix. 
	Then the approximate gradient $\g(\x)$ defined in Eq.~\eqref{eq:g} satisfies that
	\begin{equation}\label{eq:g_prop}
		\g(\x) = SS^\top \nabla \phi(\x) + S\bv \mbox{ and } \norm{\bv}^2 \le  \frac{\ell\sigma^2}{\alpha^2}, 
	\end{equation} 
	where $\bv$ is an $\ell$-dimension vector whose $i$-th entry $v^{(i)} = \frac{\zeta(\x + \alpha \si) - \zeta(\x - \alpha \si)}{2\alpha}$ and $\si$ being the $i$-th column of $S$.
\end{lemma}

Let us now demonstrate how the representation in~\eqref{eq:g_prop} motivates the need for randomized sketching matrices. 
Consider a generic gradient-based update rule:
\begin{equation}\label{eq:update}
\x_{t+1} = \x_t - \eta \cdot \g(\x_t),
\end{equation}
where $\x_t$ and $\x_{t+1}$ are the $t$-th and $(t+1)$-th iterates, $\eta > 0$ is the step size, and the estimated gradient is given by $\g(\x_t) = S_t S_t^\top \nabla \phi(\x_t) + S_t \bv_t$. 
Supposing that Assumption~\ref{ass:phi} holds, we analyze the descent behavior as follows:
\begin{align}\label{eq:motivation1}
	\phi(\x_{t+1})& 
\stackrel{\eqref{eq:update}\eqref{eq:phi}}{=}
	\phi(\x_t) - \eta \dotprod{\nabla \phi(\x_t), \g(\x_t)} + \frac{\eta^2}{2} \norm{\g(\x_t)}_A^2\\\nonumber
	\stackrel{\eqref{eq:g_prop}}{=}&
	\phi(\x_t) - \eta \dotprod{\nabla \phi(\x_t), S_tS_t^\top \nabla \phi(\x_t)} - \eta \dotprod{\nabla \phi(\x_t), S_t\bv_t} + \frac{\eta^2}{2} \norm{S_tS_t^\top \nabla \phi(\x_t) + S_t\bv_t}_A^2\\\nonumber
	\le&
	\phi(\x_t) - \eta \norm{S_t^\top \nabla \phi(\x_t)}^2 + \frac{\eta}{2} \norm{S_t^\top\nabla \phi(\x_t)}^2 + \frac{\eta}{2} \norm{\bv_t}^2 + \eta^2 \norm{S_tS_t^\top \nabla \phi(\x_t)}_A^2 + \eta^2\norm{S_t \bv_t}_A^2\\\nonumber
	=&
	\phi(\x_t) - \frac{\eta}{2} \norm{S_t^\top \nabla \phi(\x_t)}^2 + \frac{\eta}{2} \norm{\bv_t}^2 + \eta^2 \norm{S_tS_t^\top \nabla \phi(\x_t)}_A^2 + \eta^2\norm{S_t \bv_t}_A^2.
\end{align}

To ensure that $\mathbf{g}(\x_t)$ is a valid descent direction, the higher-order terms $\eta^2 \|S_t S_t^\top \nabla \phi(\x_t)\|_A^2$ and $\eta^2 \|S_t \bv_t\|_A^2$ must be properly controlled. 
Focusing on the first term as an example:
\begin{align}\label{eq:motivation2}
\|S_t S_t^\top \nabla \phi(\x_t)\|_A^2 
&= \nabla^\top \phi(\x_t) S_t S_t^\top A S_t S_t^\top \nabla \phi(\x_t) \le \|S_t^\top A S_t\|_2 \cdot \|S_t^\top \nabla \phi(\x_t)\|^2 \\\nonumber
&= \|S_t^\top A^{1/2} A^{1/2} S_t\|_2 \cdot \|S_t^\top \nabla \phi(\x_t)\|^2 = \|A^{1/2} S_t S_t^\top A^{1/2}\|_2 \cdot \|S_t^\top \nabla \phi(\x_t)\|^2.
\end{align}
The key quantity in this upper bound is $\|A^{1/2} S_t S_t^\top A^{1/2}\|_2$, which directly affects the stability and convergence behavior of the algorithm.
In the special case where $S = I_d$, corresponding to the Kiefer-Wolfowitz scheme, we require evaluations of all $d$ directional derivatives to achieve $\|A^{1/2} S_t S_t^\top A^{1/2}\|_2 = \|A\|_2$.

\textit{Thus, a central challenge is to construct a pre-defined matrix $S_t$ that can effectively control the high-order variance terms (e.g., $\eta^2 \|S_t S_t^\top \nabla \phi(\x_t)\|_A^2$) in ZOO, while requiring significantly fewer than $d$ function evaluations.}
The core insight is to interpret $A^{1/2} S_t S_t^\top A^{1/2}$ as an approximation of the product $A^{1/2} A^{1/2}$. 
This motivates the following introduction of the randomized sketching matrices, which not only offer tight spectral norm control over such matrix approximations but also significantly reduce computational costs in large-scale settings~\citep{woodruff2014sketching,avron2011randomized,cohen2016optimal}.

\subsection{Oblivious Randomized Sketching Matrix}\label{subsec:ske}

We now formally introduce the notion of an \emph{oblivious randomized sketching matrix}, which plays a central role in addressing the high-order variance challenges discussed above.
By choosing $S$ from this class of sketching matrices, we are able to tightly approximate key matrix products—such as $A^{1/2} A^{1/2}$—in spectral norm, thereby enabling control over the critical term $\|A^{1/2} S_t S_t^\top A^{1/2}\|_2$ using only a small number of directional queries.

\begin{definition}[Sketching in Matrix Product]\label{def:ske}
Given a matrix $A \in \mathbb{R}^{d \times d}$, parameters $0 < \gamma \le 1$, $k \geq 1$, and $\delta > 0$, we say a matrix $S \in \mathbb{R}^{d \times \ell}$, independent of $A$, is an oblivious $(\gamma, k, \delta)$-random sketching matrix for matrix product if the following holds with probability at least $1 - \delta$:
\begin{equation}\label{eq:mat_pd}
		\norm{A^\top S S^\top A - A^\top A}_2 \le \gamma \left(\norm{A}_2^2 + \frac{\norm{A}_F^2}{k}\right) .
	\end{equation} 
\end{definition}

According to \eqref{eq:mat_pd}, we can elaborately bound $\|A^{1/2} S_t S_t^\top A^{1/2}\|_2$ when $S_t$ is selected as an oblivious random sketching matrix. 
Specifically,
\begin{align*}
	\norm{A^{1/2} S_t S_t^\top A^{1/2} - A^{1/2}A^{1/2}}_2 \stackrel{\eqref{eq:mat_pd}}{\le} \gamma \left(\norm{A^{1/2}}_2^2 + \frac{\norm{A^{1/2}}_F^2}{k}\right).
\end{align*}
This further implies that
\begin{align}\label{eq:sta}
\norm{S_t^\top A S_t}_2 
\le
(\gamma+1)\norm{A}_2 + \gamma\frac{\tr(A)}{4k}.
\end{align}
The bound in Eq.~\eqref{eq:sta} is fundamental. 
It shows that $\|A^{1/2} S_t S_t^\top A^{1/2}\|_2$—a critical quantity in our descent analysis—can be effectively controlled using only $\ell$ directional evaluations, with $\ell \ll d$ if $\tr(A) / \norm{A}_2 \ll d$.
This provides a promising way to overcome the computational challenge discussed earlier. 
In the following, we introduce several widely used sketching matrices that satisfy Definition~\ref{def:ske}, all of which share the desirable property that the required sketching dimension $\ell$ is significantly smaller than the ambient dimension $d$.


\begin{itemize}
    \item Gaussian sketching matrix: The most classical sketching matrix is the Gaussian sketching matrix $S\in\RR^{\ell\times d}$, whose entries are i.i.d.\ from the univariate Gaussian distribution of mean 0 and variance $1/\ell$. 
Owing to the well-known concentration properties \citep{woodruff2014sketching},  Gaussian random matrices are very attractive. 
It is only requires $\ell = \Omega\left((k + \log\frac{1}{\delta}) \gamma^{-2}\right)$ to achieve the  $(\gamma, k, \delta)$-approximate matrix product property defined in Definition~\ref{def:ske} \citep{cohen2016optimal}. 

\item Rademacher sketching matrix: The Rademacher sketching matrix is also widely used in the matrix sketching whose entries are i.i.d.\ randomly sampled from $\{+1, -1\}$ and scaled with $1/\sqrt{\ell}$. 
To achieve the  $(\gamma, k, \delta)$-approximate matrix product property, it also takes $\ell = \Omega\left((k + \log\frac{1}{\delta}) \gamma^{-2}\right)$ for the Rademacher sketching matrix.
The sample size $\ell$ of the Rademacher sketching matrix is almost the same as that of the Gaussian sketching matrix.


\item Subsampled Randomized Hadamard Transform (SRHT): It is another important kind of sketching matrices. Let $H_d\in\RR^{d \times d}$ be the Walsh-Hadamard matrix with $+1$ and $-1$ entries, $D \in\RR^{d \times d}$ be a diagonal matrix with diagonal entries sampled uniformly from $\{+1,-1\}$, and $P\in\RR^{s\times d}$ be the uniform sampling matrix. 
Then when $\ell = \Omega( (k+\log(1/(\gamma \delta))\log(k/\delta))\gamma^{-2})$, the SRHT matrix $S=PH_dD$ can achieve the $(\gamma, k, \delta)$-approximate matrix product property \citep{cohen2016optimal}.

\item Sparse embedding sketching matrix: The sparse embedding matrix $S\in\RR^{\ell\times d}$ is a sparse matrix who only has $s$ 
 nonzero entries uniformly sampled from $\{+1,-1\}$ per column \citep{clarkson2017low,meng2013low,nelson2013osnap}.  For sparse subspace embeddings, the analysis in
\citet{cohen2016simpler} implies $\ell = \Omega\left(k\log\frac{k}{\delta} \gamma^{-2}\right)$ with  $s = \cO\left(\log\frac{k}{\delta} \gamma^{-1}\right)$ suffices to achieve the  $(\gamma, k, \delta)$-approximate matrix product property.
\end{itemize}

Note that our framework can support multiple sketch families is essential because they trade off \emph{generation cost} and \emph{application cost} in markedly different ways: Gaussian sketches offer strong concentration but are dense (higher arithmetic and memory traffic), Rademacher sketches are discrete and cheaper to generate and communicate, SRHT leverages fast transforms to accelerate structured multiplies, and sparse embeddings enable input-sparsity time when data or operators are sparse. These differences directly impact wall-clock time and communication overhead in large-scale deployments (e.g., LLM fine-tuning), so allowing for general sketching matrices—rather than restricting to a single Gaussian choice—delivers practical efficiency while preserving the same matrix-product approximation guarantees used in our analysis.
\section{Convergence Analysis for Convex Quadratic Function}\label{sec:convergence_convex_quad}

In this section, we begin by analyzing the behavior of ZOO methods applied to the convex quadratic objectives (see Assumption \ref{ass:phi}). 
Convex quadratic functions offer a compelling testbed for analyzing zeroth-order methods due to both their theoretical significance and practical relevance. 
From the theoretical viewpoint, it is well established that many worst-case lower bounds in convex optimization are attained by quadratic functions~\citep{nesterov2003introductory}. 
On the applied side, such functions are ubiquitous in well-established models, including classic problems like least-squares regression and its regularized variants~\citep{james2013introduction}. 
As such, studying the quadratic case yields valuable insights that generalize to broader classes of smooth convex objectives (see Section \ref{sec:convegence_smooth_convex}).
This section focuses on the classical gradient-based ZOO algorithm proposed by \citet{nesterov2017random}.
{We present it in Algorithm \ref{alg:SA} and further discuss its convergence analysis with fixed step size $\eta_t=\eta$.
}

\begin{algorithm}[t]
	\caption{Gradient-basded ZOO Algorithm~\citep{nesterov2017random}}
	\label{alg:SA}
	\begin{small}
		\begin{algorithmic}[1]
			\STATE {\bf Input:}
			Initial vector $\x_0$ and sample size $\ell$.
			\FOR {$t=0,1,2,\dots, T-1$ }
			\STATE Generate a $(1/4, k, \delta/T)$-oblivious sketching matrix $S_t \in\RR^{d \times \ell}$ for the matrix product with $\sti$ be the $i$-th column of $S_t$.
			\STATE Access to the value of $f(\cdot)$ and construct the approximate gradient
			\begin{equation*}
				\g(\x_t) = \sum_{i=1}^\ell \frac{f(\x_t+\alpha \sti) - f(\x_t - \alpha \sti)}{2\alpha} \sti.
			\end{equation*}
			\STATE Update as
			\begin{equation}
				\x_{t+1} = \x_t - \eta_t \g(\x_t),
			\end{equation}
            where $\eta_t$ is a step size at $t$.
			\ENDFOR
			\STATE {\bf Output:} $\x_T$.
		\end{algorithmic}
	\end{small}
\end{algorithm}


Note that Eqs.\eqref{eq:motivation1} and \eqref{eq:motivation2} have demonstrated the insights of how to bound the higher-order terms to improve the convergence analysis for gradient-based ZOO methods. 
This section will complete the higher-order analysis in terms of Lemma \ref{lemma:bound_phi} and \ref{lem:ske}, and justify the improved convergence results of Theorem \ref{thm:main} via the provided significant lemmas.

\begin{lemma}\label{lemma:bound_phi}
Let the objective function $\phi(\x)$ satisfy Assumption~\ref{ass:phi} and the approximate gradient $\g(\x_t)$ be defined as Eq.~\eqref{eq:g}. 
Then $\x_{t+1}$ updated as Eq.~\eqref{eq:update} satisfies that
\begin{equation} \label{eq:phi_dec}
\begin{aligned}
\phi(\x_{t+1}) 
\leq&
\phi(\x_t) - \frac{\eta}{2} \norm{S_t^\top \nabla\phi(\x_t)}^2 + \eta^2 \norm{S_t^\top A S_t}_2 \cdot \norm{S_t^\top \nabla \phi(\x_t)}^2 \\
&+ \frac{\eta}{2} \norm{\bv_t}^2 + \eta^2 \norm{S_t^\top A S_t}_2 \cdot \norm{\bv_t}^2.	
\end{aligned}
\end{equation}
\end{lemma}

\begin{lemma}\label{lem:ske}
Letting $S_t$ be a $(\frac{1}{4}, k, \delta/T)$-oblivious sketching matrix for the matrix product and matrix $A$ be positive definite, then with a probability at least $1-\delta/T$, it holds that
\begin{equation} \label{eq:sas1}
\norm{S_t^\top A S_t}_2 
\le
\frac{5\norm{A}_2}{4} + \frac{\tr(A)}{4k},\quad
\mbox{ and }\quad
\frac{1}{2} \norm{\nabla \phi(\x_t)}^2 \leq 	\norm{S_t^\top \nabla \phi(\x_t)}^2 \leq \frac{3}{2}\norm{\nabla \phi(\x_t)}^2.
\end{equation}
\end{lemma}
The above two lemmas provide a new insight by bridging the connections between bounding the higher-order term and the production approximation of a random sketching matrix.
Thus, taking full advantage of the two lemmas, we are ready to state the improved convergence analysis of Algorithm \ref{alg:SA} for the convex quadratic objective.

\begin{theorem}\label{thm:main}
Let the objective function $\phi(\x)$ satisfy Assumption~\ref{ass:phi}-\ref{ass:zeta}.
Assume that $S_t$ is a $(\frac{1}{4}, k, \frac{\delta}{T})$-oblivious sketching matrix for the matrix product.
Setting the step size $\eta_t=\eta = \left(5\norm{A}_2 + \frac{\tr(A)}{k}\right)^{-1}$,
then the sequence $\{\x_t\}$ generated by Algorithm~\ref{alg:SA} satisfies the following property with a probability at least $1 - \delta/T$,
\begin{equation}\label{eq:main}
\phi(\x_{t+1}) - \phi(\x^*) - \frac{3\ell\sigma^2}{\mu \alpha^2}
\leq
\left(1 - \frac{\mu}{4} \cdot \left(5 \norm{A}_2 + \frac{\tr(A)}{k}\right)^{-1}\right) \cdot \left(\phi(\x_t) - \phi(\x^*) - \frac{3\ell\sigma^2}{\mu \alpha^2}\right),
\end{equation}
where $\x^*$ is the optimum of the objective function $\phi(\x)$.
\end{theorem}

We now compare our theoretical results with those of \citet{yue2023zeroth} from both practical and analytical perspectives.
From an \emph{application perspective}, our work applies to a much broader class of sketching matrices. Specifically, we establish convergence guarantees for any sketch matrix satisfying the spectral approximation condition in Definition~\ref{def:ske}, including Gaussian, Rademacher, SRHT, and sparse sketches. In contrast, the analysis in \citet{yue2023zeroth} is tailored to Gaussian sketching only, and their proof techniques heavily rely on the special concentration properties of Gaussian random vectors. Consequently, their results \emph{do not directly extend} to other structured or hardware-friendly sketching matrices. 


Moreover, our results hold with \emph{high probability}, whereas \citet{yue2023zeroth} only provide convergence guarantees in \emph{expectation}. High-probability bounds are generally more informative and applicable in real-world scenarios, particularly when algorithms are deployed in noisy or mission-critical environments. To the best of our knowledge, we are the first to establish high-probability convergence guarantees for random-search ZOO algorithms with improved rates on quadratic functions~\citep{yue2023zeroth,malladi2023fine,wangcan}.

From a \emph{theoretical standpoint}, our convergence rate
$
1 - \mathcal{O}\left(\frac{\mu}{\|A\|_2 + \mathrm{tr}(A)/k}\right)
$
matches the best-known rate from \citet{yue2023zeroth} under batch size $\ell = \mathcal{O}(k)$.
However, our result improves upon theirs in several key aspects. First, our analysis features only \emph{one} perturbation term, whereas the bound in their Theorem~10 involves two sources of error:
\begin{equation}\label{eq:yue}
\begin{aligned}
&\mathbb{E}\left[\phi(\x_{t+1})\right] - \phi(\x^*) - \frac{24\tr(A)}{\mu}\left(c_1\alpha^2 + c_2\frac{\sigma^2}{\alpha^2}\right) \\
\leq& \left(1 - \frac{\mu}{24\tr(A)}\right)\left(\phi(\x_t) - \phi(\x^*) - \frac{24\tr(A)}{\mu}\left(c_1\alpha^2 + c_2\frac{\sigma^2}{\alpha^2}\right)\right),
\end{aligned}
\end{equation}
where $c_1 = \frac{5d}{16}\tr(A) + \frac{5}{384}\tr(A)$ and $c_2 = \frac{d}{3\tr(A)} + \frac{1}{72\tr(A)}$.
The presence of the extra $\mathcal{O}(\alpha^2)$ term restricts their ability to choose a large smoothing parameter $\alpha$.
In contrast, our convergence bound (see Eq.~\eqref{eq:main}) contains only the noise-induced term and no additional $\alpha^2$ penalty. This allows our method to use larger $\alpha$ to mitigate the effect of the noise term $\zeta(x)$ without incurring extra bias.
Finally, we observe that the variance term in \eqref{eq:yue} scales as
\[
\frac{24\tr(A)}{\mu} \cdot c_2 \cdot \frac{\sigma^2}{\alpha^2} = \mathcal{O}\left(\frac{d}{\mu} \cdot \frac{\sigma^2}{\alpha^2}\right),
\]
while in our result, the corresponding term scales as
$
\mathcal{O}\left(\frac{\ell}{\mu} \cdot \frac{\sigma^2}{\alpha^2}\right),
$
which is significantly smaller when $\ell \ll d$—a common case in practice. Therefore, our analysis provides both improved theoretical guarantees and greater flexibility in algorithm design and deployment.

\begin{corollary}\label{cor:main}
Suppose $\phi(\xb)$ satisfy the properties shown in Theorem~\ref{thm:main}. 
Let $S_t\in\RR^{d \times \ell}$ in Algorithm~\ref{alg:SA} be a $(\frac{1}{4}, k, \frac{\delta}{T})$- oblivious sketching matrix with $\ell = \cO\left(k + \log\frac{T}{\delta}\right)$. Assuming that $k > \log\frac{T}{\delta}$, to find a $\x_T$ such that $\phi(\x_T) - \phi(\x^*) \le \epsilon + \frac{3\ell\sigma^2}{\mu \alpha^2}$ with probability $1-\delta$, then the iteration number of Algorithm~\ref{alg:SA} requires to be 
\begin{equation}\label{eq:TT}
	T = \cO\left(\left(\frac{\norm{A}_2}{\mu} + \frac{\tr(A)}{\mu\ell }\right)\log\frac{1}{\epsilon} \right).
\end{equation}
Accordingly, the query complexity is 
\begin{equation}\label{eq:Q}
	Q = T\times\ell = \Tilde{\cO}\left(\left(\frac{\ell\norm{A}_2}{\mu} + \frac{\tr(A)}{\mu}\right)\log\frac{1}{\epsilon}\right),
\end{equation}
where $\Tilde{\cO}(\cdot)$ hides constant and $\log\log\frac{1}{\epsilon}$ terms.
\end{corollary}

\begin{remark}\label{remark_1}
{The leading term in Eq.~\eqref{eq:Q} indicates that Algorithm~\ref{alg:SA} achieves a query complexity of $\Tilde{\cO}\left(\frac{\tr(A)}{\mu}\log\frac{1}{\epsilon}\right)$, provided that the sketch size satisfies $\ell \leq \frac{\tr(A)}{\|A\|_2}$.}
For comparison, consider the approach that estimates the full gradient $\nabla \phi(\x)$ by computing all $d$ partial derivatives using zeroth-order oracle calls. Combining this with the classical convergence theory of gradient descent~\citep{nesterov2003introductory}, one requires
$
T' = \cO\left(\frac{\|A\|_2}{\mu}\log\frac{1}{\epsilon}\right)
$
iterations to reach $\epsilon$-accuracy.
Since each iteration involves $2d$ oracle calls (due to forward or central finite difference schemes), the total query complexity becomes
\[
Q' = T' \times 2d = \cO\left( \frac{d \|A\|_2}{\mu} \log\frac{1}{\epsilon} \right).
\]
When the eigenvalues of $A$ are highly non-uniform—i.e., $A$ has a skewed spectrum—it is often the case that
\[
\frac{\tr(A)}{\mu} \ll \frac{d \|A\|_2}{\mu}.
\]
This phenomenon has been shown in several studies to demonstrate the non-uniform property in the least-squares problems~\citep{deng2012mnist} and deep models~\citep{sagun2016eigenvalues,yao2020pyhessian}.
Therefore, our algorithm offers substantial query efficiency over the classical zeroth-order gradient approximation method in such scenarios, while still achieving comparable convergence guarantees.
Moreover, Eq.~\eqref{eq:TT} reveals that increasing the sketch size $\ell$ can reduce the iteration complexity of our algorithm, up to the threshold $\ell = \frac{\tr(A)}{\|A\|_2}$. Beyond this point, however, Eq.~\eqref{eq:Q} shows that the total query complexity no longer decreases with larger $\ell$.
In fact, if $\ell > \frac{\tr(A)}{\|A\|_2}$, then further increasing the sketch size will lead to a higher query complexity, as each iteration incurs more function evaluations without yielding proportional gains in convergence speed.
\end{remark}

\begin{remark}
We emphasize that Corollary~\ref{cor:main} is established specifically for the Gaussian sketching matrix and the Rademacher sketching matrix.
The key reason for this restriction lies in the fact that both sketching matrices satisfy strong concentration properties, which allow to construct a $(\frac{1}{4}, k, \frac{\delta}{T})$- oblivious sketching matrix with $\ell = \cO\left(k + \log\frac{T}{\delta}\right)$.
For other sketching matrices—such as the SRHT or sparse embedding matrices—similar bounds can be derived for different sampling size orders of $\ell$.
We refer interested readers to pursue such extensions based on the general sketching framework introduced in Section~\ref{subsec:ske}.
\end{remark}

\section{Hessian-Aware ZOO with Randomized Sketching}\label{sec:hessian_aware}

A central theme in ZOO is to leverage curvature information to accelerate convergence without explicit gradients or Hessians; see, e.g., Hessian-aware evolutionary strategies and natural-gradient methods~\citep{ye2018hessian,wierstra2014natural,hansen2016cma}. 
However, as highlighted in our introduction (Open Questions of \citet{scheinberg2022finite}), directly combining quasi-Newton ideas with RFD gradients is problematic due to estimator inaccuracy and the ensuing instability of line search and curvature updates. 
In this section, we show that \emph{oblivious randomized sketching} provides a principled remedy: by embedding curvature through a sketched preconditioner and controlling matrix products in spectral norm, we obtain a Hessian-aware ZOO framework that (i) retains the low per-iteration query cost of randomized search, (ii) admits high-probability descent guarantees, and (iii) achieves query complexity scaling with $\Tilde{\cO}\!\left(\max_t\bigl\{\rho_t^{-1}\,\tr(H_t^{-1}A)\bigr\}\,\log \tfrac{1}{\epsilon}\right)$, thereby addressing both the variance and dimension-dependence barriers.

In this section, we consider the case that an approximated Hessian matrix $H_t$ can be obtained.  
We will study how the approximated Hessian matrix $H_t$ affects the convergence rate when the curvature information embedded in $H_t$  is exploited in the zeroth-order algorithm.
To exploit the curvature information embedded in the approximate Hessian, we design the descent direction as 
\begin{equation}\label{eq:tg}
	\tg(\x_t) = \frac{1}{\ell}\sum_{i=1}^\ell \frac{f(\x_t + \alpha H_t^{-1/2}\sti) - f(\x_t - \alpha H_t^{-1/2}\sti)}{2\alpha} H_t^{-1/2} \sti.
\end{equation} 
The descent direction $\tg(\x_t)$ in above equation has been used in several works \citep{nesterov2017random,ye2018hessian}.
Given the descent direction $\tg(\x_t)$, we can use it to update the sequence $\{\x_t\}$.
The detailed description is listed in Algorithm~\ref{alg:alg2}.

\begin{algorithm}[t]
	\caption{Hessian-aware Zeroth-order Algorithm\citep{ye2018hessian}}
	\label{alg:alg2}
	\begin{small}
		\begin{algorithmic}[1]
			\STATE {\bf Input:}
			Initial vector $x_0$, sample size $\ell$, and step size sequence $\{\eta_t\}$ and Hessian approximation sequence $\{H_t\}$.
			\FOR {$t=0,1,2,\dots, T-1$ }
			\STATE Generate a $(1/4, k, \delta/T)$-oblivious sketching matrix $S_t \in\RR^{d \times \ell}$ for the matrix product with $\sti$ be the $i$-th column of $S_t$.
			\STATE Access to the value of $f(\cdot)$ and construct the approximate gradient
			\begin{equation*}
				\tg(\x_t) =  \frac{1}{\ell}\sum_{i=1}^\ell \frac{f(\x_t + \alpha H_t^{-1/2}\sti) - f(\x_t - \alpha H_t^{-1/2}\sti)}{2\alpha} H_t^{-1/2} \sti.
			\end{equation*}
			\STATE Update as
			\begin{equation}\label{eq:update1}
				\x_{t+1} = \x_t - \eta_t \tg(\x_t).
			\end{equation}
			\ENDFOR
			\STATE {\bf Output:} $\x_T$.
		\end{algorithmic}
	\end{small}
\end{algorithm}

Similar to Lemma~\ref{lem:g_prop}, we can represent $\tg(\x_t)$ in a novel form.
\begin{lemma}\label{lem:tg_prop}
	Let the objective function $\phi(\x)$ satisfy Assumption~\ref{ass:phi} and Assumption~\ref{ass:zeta}. 
	$S_t\in\RR^{d \times \ell}$ is a pre-defined matrix.
	Then the descent direction $\tg(\x_t)$ defined in Eq.~\eqref{eq:tg} satisfies that
	\begin{equation}\label{eq:tg_prop}
		\tg(\x_t) = H_t^{-1/2}S_tS_t^\top  H_t^{-1/2}\nabla\phi(\x_t) + H_t^{-1/2} S_t\tv_t \mbox{ and } \norm{\tv_t}^2 \le  \frac{\ell\sigma^2}{\alpha^2}, 
	\end{equation} 
	where $\tv_t$ is a $\ell$-dimension vector whose $i$-th entry $\tv_t^{(i)} = \frac{\zeta(\x_t + \alpha H_t^{-1/2} \sti) - \zeta(\x_t - \alpha H_t^{-1/2} \sti)}{2\alpha}$ and $\sti$ being the $i$-th column of $S_t$.
\end{lemma}

If the matrix $S_tS_t^\top$ provides a good approximation of the identity matrix, Eq.~\eqref{eq:tg_prop} indicates that $\tg(\x_t)$ can well approximate $H_t^{-1} \nabla \phi(\x_t)$, which corresponds to the exact descent direction in the approximate Newton method with $H_t$ as the approximate Hessian \citep{ye2021approximate}. Hence, Algorithm~\ref{alg:alg2} can be interpreted as a second-order algorithm that only requires zeroth-order oracle access.
Moreover, we note that Eq.~\eqref{eq:g_prop} arises as a special case of Eq.~\eqref{eq:tg_prop} when $H_t$ is the identity matrix. This observation suggests that Algorithm~\ref{alg:alg2} is expected to exhibit a convergence behavior similar to that of Algorithm~\ref{alg:SA}, whose iteration complexity depends on the weakly dimension-dependent factor $\tr(A) / \mu$, rather than the standard $d\cdot\frac{L}{\mu}$.
Thus, in the following theorem, we show Algorithm~\ref{alg:alg2} has a convergence rate between 
one of the zeroth-order approximate Newton methods and one of Algorithm~\ref{alg:SA}.

\begin{theorem}
\label{thm:main1}
	Let the objective function $\phi(\x)$ satisfy Assumption~\ref{ass:phi}-\ref{ass:zeta}.
	Assume that $S_t$ is a $(\frac{1}{4}, k, \frac{\delta}{T})$-oblivious sketching matrix for the matrix product, and the approximate Hessian $H_t$ satisfies that
	\begin{equation}\label{eq:preceq}
		\rho_t H_t \preceq A .
	\end{equation}
	Setting the step size $\eta_t = \left(  5\norm{ H_t^{-1/2} A H_t^{-1/2} }_2 + \frac{\tr(H_t^{-1} A)}{k}\right)^{-1}$, then the sequence $\{\x_t\}$ generated by Algorithm~\ref{alg:alg2} satisfies the following property holding with a probability at least $1 - \delta/T$,
 \begin{small}
\begin{equation}\label{eq:main1}
	\phi(\x_{t+1}) -  \phi(\x^*) - \frac{3\ell\sigma^2}{\rho_t \alpha^2} 
	\leq
	\left(1 - \frac{\rho_t}{4} \left( \norm{ H_t^{-1/2} A H_t^{-1/2} }_2 + \frac{\tr(H_t^{-1} A)}{k} \right)^{-1} \right)\cdot \left(\phi(\x_t) - \phi(\x^*) - \frac{3\ell\sigma^2}{\rho_t \alpha^2}\right).
\end{equation}    
\end{small}
\end{theorem}

\begin{corollary}\label{cor:main1}
Suppose that $\phi(\x)$ satisfies the properties shown in Theorem \ref{thm:main1}.
	Let $S_t\in\RR^{d \times \ell}$ in Algorithm~\ref{alg:alg2} be a $(\frac{1}{4},k,\frac{\delta}{T})$-oblivious sketching matrix with $\ell = \cO\left(k + \log\frac{T}{\delta}\right)$. Assuming that $k > \log\frac{T}{\delta}$ and denoting $\rho = \min_t\rho_t\neq 0$, to find a $\x_T$ such that $\phi(\x_T) - \phi(\x^*) \le \epsilon + \frac{3\ell\sigma^2}{\rho \alpha^2}$ with probability $1-\delta$, then the iteration number of Algorithm~\ref{alg:SA} requires to be 
	\begin{equation}\label{eq:TTT}
		T = \cO\left(\max_t\left\{\rho_t^{-1}\left( \norm{ H_t^{-1/2} A H_t^{-1/2} }_2 + \frac{\tr(H_t^{-1} A)}{k} \right)\right\}\log\frac{1}{\epsilon} \right).
	\end{equation}
	Accordingly, the query complexity is 
	\begin{equation}\label{eq:QQ}
		Q = \Tilde{\cO}\left(\max_t\left\{\rho_t^{-1}\left( \norm{ H_t^{-1/2} A H_t^{-1/2} }_2 \cdot \ell + \tr(H_t^{-1} A) \right)\right\}\log\frac{1}{\epsilon} \right).
	\end{equation}
\end{corollary}

Incorporating curvature in ZOO has a long history, ranging from natural-gradient and evolution-strategy methods to recent randomized second-order schemes; see, e.g., \citet{wierstra2014natural,hansen2016cma,ye2018hessian,gower2019rsn,cartis2023scalable,ye2023mirror}. These approaches share the ambition of leveraging (explicit or implicit) Hessian information to accelerate convergence while relying only on function evaluations.

    On the positive side, Hessian-aware designs typically achieve improved iteration complexity by adapting steps to local curvature.
However, many existing analyses either (i) rely on full or blockwise Hessian information whose construction costs at least $\cO(\ell^2)$ function queries per iteration (due to forming terms like $S_t^\top \nabla^2\phi(\x_t) S_t$), or (ii) establish guarantees only in expectation and for specific sketch distributions (most often Gaussian), which complicates extensions to structured or communication-efficient sketches.

When the batch size satisfies
\[
\ell \;\le\; \frac{\tr(H_t^{-1}A)}{\|H_t^{-1/2} A H_t^{-1/2}\|_2},
\]
our Algorithm~\ref{alg:alg2} attains the high-probability query complexity
\[
\Tilde{\cO}\!\left(\max_t\bigl\{\rho_t^{-1}\,\tr(H_t^{-1}A)\bigr\}\,\log \tfrac{1}{\epsilon}\right),
\]
as shown by Eq.~\eqref{eq:QQ}. 
This bound scales with the \emph{trace} surrogate $\tr(H_t^{-1}A)$ and maintains per-iteration query cost $\cO(\ell)$, since we avoid forming second-order blocks explicitly. Let us compare our main results with other ZOO methods incorporating curvature information.

\begin{itemize}
   \item \emph{Comparison with \citet{ye2018hessian}.}
The theory of \citet{ye2018hessian} yields a query complexity of
\[
\mathcal{O}\!\left(d\cdot \max_t\bigl\{\rho_t^{-1}\,\|H_t^{-1/2} A H_t^{-1/2}\|_2\bigr\}\cdot \log \tfrac{1}{\epsilon}\right),
\]
which replaces the trace term by $d$ times the spectral term and thus can be significantly larger when the spectrum of $H_t^{-1/2} A H_t^{-1/2}$ is skewed. 
In regimes where
\[
\rho_t^{-1}\,\tr(H_t^{-1} A)\;\ll\; d\cdot \rho_t^{-1}\, \|H_t^{-1/2} A H_t^{-1/2}\|_2,
\]
our result provides a much tighter query complexity.

\item \emph{Comparison with randomized subspace Newton~\citep{gower2019rsn}.}
This method computes the search direction
\[
d_t \;=\; S_t\!\left[S_t^\top \nabla^2 \phi(\x_t) S_t \right]^{-1}\! S_t^\top \nabla \phi(\x_t),
\]
which requires $\cO(\ell)$ function queries to approximate $S_t^\top \nabla \phi(\x_t)$ but \emph{$\cO(\ell^2)$} queries to approximate the projected Hessian $S_t^\top \nabla^2 \phi(\x_t) S_t$ per iteration. 
While a linear rate is established, the per-iteration query cost grows quadratically in $\ell$, whereas our sketched framework keeps it linear in $\ell$ by controlling matrix products in spectral norm without explicitly forming second-order blocks.

\item \emph{Comparison with \citet{cartis2023scalable}.}
For nonlinear least-squares problem $\phi(\x)=\sum_{i=1}^m r_i(\x)^2$, one can approximate
\[
S_t^\top \nabla^2 \phi(\x_t) S_t \;\approx\; (J(\x_t) S_t)^\top (J(\x_t) S_t),
\]
with $J(\x_t)$ the Jacobian of residuals.
When each $r_i(\x)$ (and hence $J$) is accessible, \citet{cartis2023scalable} obtain favorable query efficiency. 
However, these assumptions are specific to least-squares structure and do not generally hold for black-box objectives, limiting applicability beyond that class.

\item \emph{Comparison with iterative Hessian approximation~\citep{hansen2016cma,ye2023mirror}.}
CMA-ES maintains and updates a covariance/Hessian surrogate using only function evaluations and is conjectured to converge to the true Hessian empirically~\citep{hansen2016cma}. 
\citet{ye2023mirror} provide an iterative Hessian-approximation scheme with provable convergence to the true Hessian under suitable conditions.
These techniques avoid explicit second-order blocks at each step, but existing analyses typically do not furnish high-probability {query-complexity} guarantees comparable to our trace-controlled bounds, nor do they directly leverage oblivious sketching to obtain dimension-robust rates with structured sketches.
\end{itemize}

Across these Hessian-aware lines of work, our contribution is to show that \emph{oblivious randomized sketching} yields:
(i) high-probability convergence with query complexity scaling as $\Tilde{\cO}\!\bigl(\max_t\{\rho_t^{-1}\tr(H_t^{-1}A)\}\log(1/\epsilon)\bigr)$ under mild sketch-size conditions;
(ii) linear-in-$\ell$ per-iteration query cost (no $\ell^2$ Hessian blocks);
and (iii) applicability to a broad family of sketches (Gaussian, Rademacher, SRHT, sparse), enabling communication- and hardware-friendly implementations beyond Gaussian-only analyses.

\begin{remark}
	Now, we will discuss how the precision of the approximate Hessian $H_t$ affects the query complexity of Algorithm~\ref{alg:alg2}.
	First, we consider the case that $H_t$ well approximates the Hessian matrix $A$. 
	For example, if it holds that $\frac{1}{2}H_t \preceq A\preceq \frac{3}{2} H_t$, Eq.~\eqref{eq:main1} shows that Algorithm~\ref{alg:alg2} converges with a rate $\cO(1 - d^{-1})$ which is independent of the condition number of the objective function.
	Accordingly, the query complexity of Eq.~\eqref{eq:QQ} becomes $\Tilde{\cO}\left(d\log\frac{1}{\epsilon}\right)$ which depends on the dimension $d$ but without the condition number $\frac{L}{\mu}$. 
	
	On the other hand, if $H_t$ approximates the Hessian matrix $A$ poorly, for an extreme case that $H_t$ equals the identity matrix, then $\rho_t$ will be $\mu$. 
	The query complexity in Eq.~\eqref{eq:QQ} will reduce to the one in Eq.~\eqref{eq:Q} which depends on the condition number $L/\mu$ and a weak dimension $\tr(A)/ L$ just as we represent Eq.~\eqref{eq:Q} as $\Tilde{\cO}\left(\frac{L}{\mu}\cdot \frac{\tr(A)}{L}\log\frac{1}{\epsilon}\right)$.
	
	In summary, the approximate Hessian $H_t$ can help to reduce the dependency of the condition number.
	However, this maybe is at the cost of increasing the dependency of one the dimension. 
	Our theory shows Algorithm~\ref{alg:alg2} may achieve a balance between the dimension and the condition number if the approximate Hessian $H_t$ is properly constructed.
\end{remark}

\section{Convergence Analysis for Smooth and Strongly Convex Function}
\label{sec:convegence_smooth_convex}

In the previous sections, we demonstrated that our sketch-based framework enables ZOO for convex quadratic objectives, with convergence guarantees depending only on $\tr(A)/\mu$ rather than the ambient dimension $d$. This result reveals the potential of reducing query complexity through sketching when the problem structure is favorable. 
To broaden the applicability of our method, we now extend our analysis to general $L$-smooth and $\mu$-strongly convex functions whose Hessians are $\mL$-Lipschitz continuous. Specifically, we show that when the curvature of the objective function, as quantified by $\tr(\nabla^2 \phi(\x))/\mu$, is significantly smaller than $d\cdot \frac{L}{\mu}$, our Algorithm~\ref{alg:SA} can achieve a strictly better query complexity than the classical zeroth-order methods, which typically require $\cO\left(d\frac{L}{\mu}\log\frac{1}{\epsilon}\right)$ function evaluations.

Let us give the representation of the approximation gradient $g(\x_t)$ first.

\begin{lemma}\label{lem:g_prop1}
	Let the objective function $\phi(\x)$ satisfy Assumption~\ref{ass:mu}-\ref{ass:mL}. 
	$S_t\in\RR^{d \times \ell}$ is a pre-defined matrix.
	Then the approximate gradient $g(\x_t)$ defined in Eq.~\eqref{eq:g} satisfies that
	\begin{equation}\label{eq:g_prop1}
		g(\x_t) = S_tS_t^\top \nabla \phi(\x_t) +S_t \bw_t + S_t\bv_t, \quad \norm{\bw_t}^2 \le \frac{\mL^2\alpha^4}{36} \norm{S_t}_F^6, \mbox{ and } \norm{\bv_t}^2 \le  \frac{\ell\sigma^2}{\alpha^2}, 
	\end{equation} 
	where $\bw_t$ and $\bv_t$ are  $\ell$-dimension vectors whose $i$-th entries are $\bw_t^{(i)} = \frac{\phi(\x_t+\alpha \sti) - \phi(\x_t-\alpha \sti) - 2 \alpha \langle\nabla \phi(\x_t), \sti\rangle}{2 \alpha}$ and $\bv_t^{(i)} = \frac{\zeta(\x_t + \alpha \sti) - \zeta(\x_t - \alpha \sti)}{2\alpha}$ respectively. And $\sti$ is the $i$-th column of $S_t$.
\end{lemma}

We provide the convergence analysis of Algorithm \ref{alg:SA} for the smooth and strongly convex objective function.

\begin{lemma}\label{lemma:th_part1}
	Let the objective function $\phi(\x)$ satisfy Assumption~\ref{ass:mu}-\ref{ass:mL}. 
We assume that  $S_t\in\RR^{d \times \ell}$ is an oblivious $(1/4, k, \delta/T)$-sketching matrix for the matrix product. 
By setting the step size $\eta_t=\eta = \left(5+\frac{d}{k}\right)^{-1} L^{-1}$, then the sequence $\{\x_t\}$ generated by Algorithm~\ref{alg:SA} satisfies the following property holding with a probability at least $1 - \delta/T$,
\begin{small}
\begin{equation}\label{eq:dec}
\phi(\x_{t+1}) -\phi(\x^*)
\leq 
\left(1 - \frac{\mu}{4L} \left(5 + \frac{d}{k}\right)^{-1}\right) \cdot \Big( \phi(\x_t) - \phi(\x^*) \Big) +  \frac{3}{2L}\left(5 + \frac{d}{k}\right)^{-1} \cdot \left(\frac{\mL^2\alpha^4}{36} \norm{S_t}_F^6+\frac{\ell\sigma^2}{\alpha^2}\right).
\end{equation}
\end{small}
\end{lemma}

\begin{lemma}\label{lemma:th_part2}
Let the objective function $\phi(\x)$ satisfy Assumption~\ref{ass:mu}-\ref{ass:mL}.
Assume that $S_t$ is a $(\frac{1}{4}, k, \delta/T)$-oblivious sketching matrix for the matrix product.
Let us denote $C_1 = \max_{\x} 5\norm{\nabla^2 \phi(\x)}_2 + \frac{\tr(\nabla^2\phi(\x))}{k}$ and $C_2 = \sqrt{\frac{5}{4} + \frac{d}{4k}} $. 
If $\x_t$ is in the local range that $(\phi(\x_t) - \phi(\x^*))^{1/2} \le \frac{3\sqrt{2}C_1^2}{16\sqrt{L}\mL C_2^3}$, 
setting the step size $\eta = 1/(2C_1)$,
then the sequence $\{\x_t\}$ generated by Algorithm~\ref{alg:SA} satisfies the following property holding with a probability at least $1 - \delta/T$,
\begin{equation}\label{eq:dec1}
	\phi(\x_{t+1}) - \phi(\x^*) - C_4 \le \left(1 - \frac{\mu}{8C_1}\right)\cdot \Big(\phi(\x_t) - \phi(\x^*) - C_4\Big),
\end{equation}
where $C_3$ is defined in Lemma \ref{lem:bound_S_F_norm} and $C_4 = \mu^{-1}\left(\frac{20\ell\sigma^2}{\alpha^2} + \frac{5C_3^3\mL^2\alpha^4}{9}\right)+ \frac{8C_2^3\mL}{3C_1^2\mu}\left( \frac{C_3^{9/2}\mL^3\alpha^6}{6^3} + \frac{\ell^{3/2}\sigma^3}{\alpha^3} \right)$.

\end{lemma}

Combining the above two lemmas, we can obtain the following main theorem.

\begin{theorem}\label{thm:main2}
Let the objective function $\phi(\x)$ satisfy Assumption~\ref{ass:mu}-\ref{ass:mL}. 
We assume that  $S_t\in\RR^{d \times \ell}$ is an oblivious $(1/4, k, \delta/T)$-sketching matrix for the matrix product. 
If $\x_t$ satisfies $(\phi(\x_t) - \phi(\x^*))^{1/2} > \frac{3\sqrt{2}C_1^2}{16\sqrt{L}\mL C_2^3}$, setting the step size $\eta_t=\eta = \left(5 + \frac{d}{k}\right)^{-1} L^{-1}$, then Algorithm~\ref{alg:SA} has a convergence property shown in Eq.~\eqref{eq:dec} with a probability at least $1 - \delta/T$. 
If $\x_t$ enters into the local region that satisfies $(\phi(\x_t) - \phi(\x^*))^{1/2} \leq \frac{3\sqrt{2}C_1^2}{16\sqrt{L}\mL C_2^3}$, setting the step size $\eta_t=\eta = 1/(2C_1)$,
then Algorithm~\ref{alg:SA} has a convergence property shown in Eq.~\eqref{eq:dec1} with a probability at least $1 - \delta/T$,. 
\end{theorem}

\begin{corollary}\label{cor:main2}
Let the objective function satisfy the properties required in Theorem~\ref{thm:main2}. 
	Let $S_t\in\RR^{d \times \ell}$ in Algorithm~\ref{alg:alg2} be a $(\frac{1}{4},k,\frac{\delta}{T})$-oblivious sketching matrix with $\ell = \cO\left(k + \log\frac{T}{\delta}\right)$. Assuming that $k > \log\frac{T}{\delta}$, to find a $\x_T$ such that $\phi(\x_T) - \phi(\x^*) \le \epsilon + C_4$ with a probability at least $1 - \delta$, then the iteration number of Algorithm~\ref{alg:SA} requires to be 
	\begin{equation}\label{eq:T3}
		T = \cO\left(\frac{d}{\ell}\cdot \frac{L}{\mu} + \frac{C_1}{\mu}\log\frac{1}{\epsilon}\right).
	\end{equation}
	Accordingly, the query complexity is 
	\begin{equation}\label{eq:Q3}
		Q = \Tilde{\cO}\left(\mu^{-1}dL +  \mu^{-1} \left( \max_\x \left(\norm{\nabla^2\phi(\x)}_2\cdot \ell + \tr(\nabla^2\phi(\x))\right)\log\frac{1}{\epsilon}\right) \right).
	\end{equation}
\end{corollary}

\begin{remark}
The communication complexity in Eq.~\eqref{eq:Q3} has two terms. 
The first one is $\mu^{-1}dL$ which is independent of the target precision $\epsilon$. 
This term can be regarded as an overhead before reaching the local region.
If the sample size $\ell\leq \frac{\tr(\nabla^2\phi(\x))}{\norm{\nabla^2 \phi(\x)}_2}$, the
the second term will reduce to  $\Tilde{\cO}\left(\max_\x\frac{\tr(\nabla^2(\phi(\x)))}{\mu}\cdot\log\frac{1}{\epsilon}\right)$ which dominates the query complexity.
Thus, Eq.~\eqref{eq:Q3} shows that if $\tr(\nabla^2\phi(\x)) \ll dL$, then Algorithm~\ref{alg:SA} with properly chosen step sizes can achieve a lower query complexity than $\cO\left(\frac{dL}{\mu}\log\frac{1}{\epsilon}\right)$.
Finally, we should emphasize that the results in Theorem \ref{thm:main2} could be generalized to the Hessian-aware version by the similar techniques in Section \ref{sec:hessian_aware}. We leave these discussions for the interested readers.   
\end{remark}

\section{Trace Estimation via Zeroth-Order Oracle}\label{sec:trace_estimation}

In the previous sections, we have established that achieving a fast convergence rate for ZOO requires selecting a step size on the order of $\left(\tr(\nabla^2 \phi(\x))/\ell\right)^{-1}$. However, since direct access to the Hessian $\nabla^2 \phi(\x)$ is unavailable in the zeroth-order setting, we now present a method to estimate its trace using only function value queries.

Our approach is inspired by the second-order expansion of $\phi(\x)$:
\begin{align*}
	\phi(\x_t+\alpha \sti) =& \phi(\x_t) + \alpha\dotprod{\nabla \phi(\x_t), \sti} + \frac{\alpha^2}{2} (\sti)^\top A \sti, \\
	\phi(\x_t-\alpha \sti) =&\phi(\x_t) - \alpha\dotprod{\nabla \phi(\x_t), \sti} + \frac{\alpha^2}{2} (\sti)^\top A \sti.
\end{align*}
where $A$ denotes the (constant) Hessian matrix in the quadratic case. By averaging the two expressions, the first-order terms cancel, yielding:
\begin{equation}
	\frac{\phi(\x_t+\alpha \s_i) + \phi(\x_t-\alpha \s_i) - 2\phi(\x_t)}{\alpha^2} = \s_i^\top A \s_i.
\end{equation}
This identity provides a mechanism to estimate individual quadratic forms involving $A$, and by leveraging classical randomized trace estimation techniques in Lemma \ref{lem:trace} ~\citep{avron2011randomized,cortinovis2021randomized,roosta2015improved,boutsidis2013improved}, the trace of $A$ can be approximated as:
\begin{equation*}
	\tr(A) \approx \frac{1}{\alpha^2} \sum_{i=1}^\ell \left[ \phi(\x_t+\alpha \s_i) + \phi(\x_t-\alpha \s_i) - 2\phi(\x_t) \right]=\sum_{i=1}^{\ell}(\sti)^\top A \sti.
\end{equation*}

While the derivation above assumes $\phi(\x)$ is quadratic, our method extends to more general settings. Specifically, we consider convex objective functions whose Hessians are $\mL$-Lipschitz continuous. This extension allows us to replace $\phi(\x)$ with the accessible objective $f(\x)$ in the estimation formula, thereby yielding a practical trace estimate based solely on zeroth-order access. The formal guarantee for this estimate is presented in the following theorem.

\begin{theorem}\label{thm:trace}
Let $S_t\in\RR^{d\times\ell}$ be the Gaussian, Rademacher, or SRHT sketching matrix, $0<\delta<1$ and $0<\varepsilon<1$. Assume that $\phi(\x)$ is convex and twice differentiable. Let Assumption~\ref{ass:zeta} and~\ref{ass:mL} hold. Define
\begin{equation}\label{eq:H}
	\tau(\x_t, S_t) = \sum_{i=1}^{\ell} \frac{f(\x_t+\alpha \sti) + f(\x_t-\alpha \sti) - 2f(\x_t)}{ \alpha^2}.
\end{equation}
Then,
\begin{enumerate}
	\item if $S_t$ is the Gaussian sketching matrix and $\ell = 8\varepsilon^{-2}\frac{\norm{\nabla^2\phi(x)}_2}{\tr(\nabla^2\phi(x))}\log\frac{2}{\delta}$,
	\item  if $S_t$ is the Rademacher sketching matrix and $\ell = 8\varepsilon^{-2}(1+\varepsilon)\frac{\norm{\nabla^2\phi(x)}_2}{\tr(\nabla^2\phi(x))}\log\frac{2}{\delta}$,
	\item if $S_t$ is the SRHT matrix and $\ell = \left(1 + \sqrt{8\ln\frac{2n}{\delta}}\right)^2\varepsilon^{-2}\log\frac{2}{\delta}$,	
\end{enumerate}
it holds that with a probability at least $1-\delta$,
\begin{equation}\label{eq:H_tr}
	\Big| \tau(\x_t, S_t) - \tr\left(\nabla^2 \phi(\x_t)\right) \Big| \le \varepsilon \cdot \tr\left(\nabla^2 \phi(\x_t)\right) + \frac{\mL\sum_{i=1}^\ell\norm{\sti}^3}{3}\alpha + \frac{4\ell\sigma}{\alpha^2},
\end{equation}
with 
\begin{equation*}
	\sum_{i=1}^\ell\norm{\sti}^3 
	\le
	\begin{cases}
		d^{3/2}\ell^{-1/2} & \mbox{ if $S_t$ is the  Rademacher or SRHT sketching matrix, }\\
		\left(2d + \log\frac{\ell}{\delta}\right)^{3/2} \ell^{-1/2} & \mbox{ if $S_t$ is the Gaussian sketching matrix. }
	\end{cases} 
\end{equation*}
\end{theorem}

Recall that Eq.~\eqref{eq:g} defines the approximation of the gradient via randomized sketching. 
Building on this idea, Eq.~\eqref{eq:H} further extends the same sketching framework to approximate the Hessian information. 
Comparing Eq.~\eqref{eq:H} with Eq.~\eqref{eq:g}, we observe that the same set of queries used for gradient estimation can also be reused to estimate the trace of the Hessian matrix $\nabla^2 \phi(\x)$. 
As a result, the Hessian trace can be obtained at little additional cost in terms of function value queries.  

Theorem~\ref{thm:trace} establishes that a sketching size of 
\[
\ell = \cO\!\left(\frac{\norm{\nabla^2\phi(\x)}_2}{\tr(\nabla^2\phi(\x))}\log\frac{1}{\delta}\right) 
     = \cO\!\left(\log\frac{1}{\delta}\right)
\] 
is sufficient to guarantee Eq.~\eqref{eq:H_tr} with $\varepsilon = 1/2$, when the sketching matrix is Gaussian or Rademacher. 
This implies that  
\begin{align*}
\frac{1}{2}\,\tr\!\left(\nabla^2 \phi(\x)\right) 
\;\lesssim\; \tau(\x, S_t) 
\;\lesssim\; \frac{3}{2}\,\tr\!\left(\nabla^2 \phi(\x)\right).
\end{align*}
Consequently, $\tau(\x, S_t)$ can be employed to determine the step size in Algorithm~\ref{alg:SA}, thereby ensuring the convergence rate stated in Eq.~\eqref{eq:main}.  

Although Theorem~\ref{thm:trace} specifically addresses Gaussian, Rademacher, and SRHT sketching matrices, similar guarantees extend to a broader class of randomized sketching matrices with sub-Gaussian properties, as implied by Lemma~\ref{lem:trace1}. 
Since the sketching matrices defined in Definition~\ref{def:ske} are typically sub-Gaussian, they inherit the approximation properties highlighted in Theorem~\ref{thm:trace}.  

\section{Experiments}
\label{sec:experiment}

In this section, we conduct experiments to validate the theoretical results established in the preceding sections empirically.
We begin by evaluating Algorithm~\ref{alg:SA} using various choices of sketching matrices on a synthetic convex quadratic problem. This setting allows us to directly assess the convergence behavior and dimension dependence according to our theory.
Next, we test our method on logistic regression tasks using real-world datasets, which serve to demonstrate the effectiveness of our approach under general $L$-smooth and $\mu$-strongly convex functions as discussed in Section~\ref{sec:convegence_smooth_convex}.
For simplicity and to isolate the core behavior of our algorithm, we set $\zeta(\x) = 0$ in all experiments. This simplification does not affect the generality of our theoretical claims or the validity of the experimental results.

\subsection{Synthetic Convex Quadratic Optimization}

We first evaluate our method on synthetic quadratic functions where $\phi(\x)$ is defined as
\begin{align*}
	\phi(\x) = \frac{1}{2}\x^\top A \x + \frac{\lambda}{2}\norm{\x}^2 - \dotprod{\x, \mathbf{a}},
\end{align*}
with $A \in \RR^{d \times d}$ being a positive definite matrix, $d = 300$, $\lambda = 10^{-4}$ a scalar parameter, and $\mathbf{a} \in \RR^d$ a given vector. 

In our experiments, we construct $A = U \Lambda U^\top$, where $U$ is a $d \times d$ orthonormal matrix and $\Lambda$ is a diagonal matrix whose entries $\lambda_{i} > 0$ represent the eigenvalues of $A$. We consider the following three types of eigenvalue distributions for $\Lambda$:

\begin{itemize}
    \item \textit{Exponential decay:} The eigenvalues decay exponentially with rate $0.95$. In this case, we have $\frac{\tr(A) + d\lambda}{L} \approx 20$.
    
    \item \textit{Polynomial decay ($1/i$):} The $i$-th eigenvalue is set to $1/i$. Here, $\frac{\tr(A) + d\lambda}{L} \approx 6$.
    
    \item \textit{Polynomial decay ($1/\sqrt{i}$):} The $i$-th eigenvalue is set to $1/\sqrt{i}$, resulting in $\frac{\tr(A) + d\lambda}{L} \approx 33$.
\end{itemize}

\begin{figure*}[!ht]
	\begin{center}
		\centering
		\subfigure[\textsf{Eigenvalue distribution for ``exponential decay''.}]{\includegraphics[width=55mm]{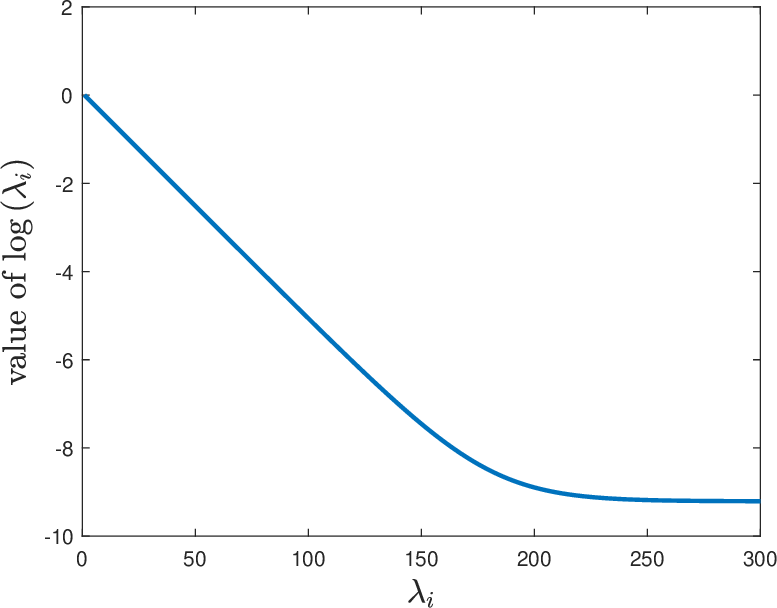}}~
		\subfigure[\textsf{Eigenvalue distribution for ``poly decay''.}]{\includegraphics[width=55mm]{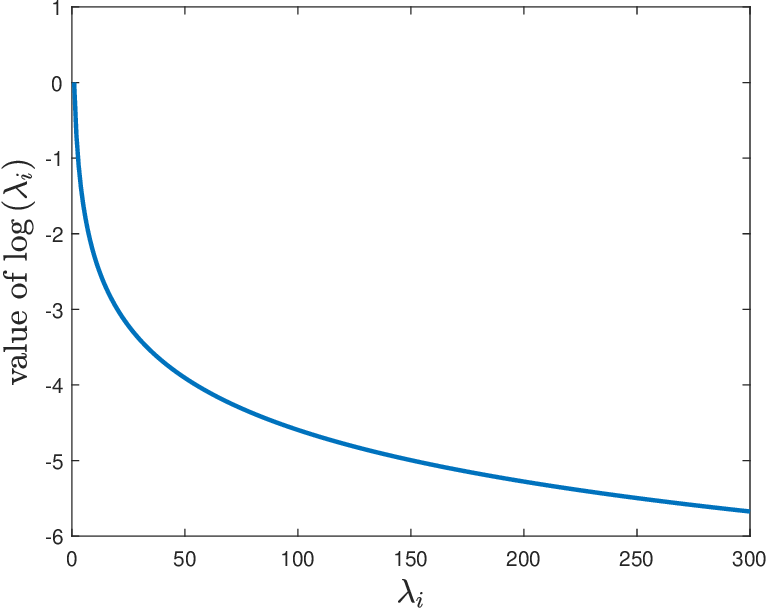}}~
		\subfigure[\textsf{Eigenvalue distribution ``poly decay sqrt''}]{\includegraphics[width=55mm]{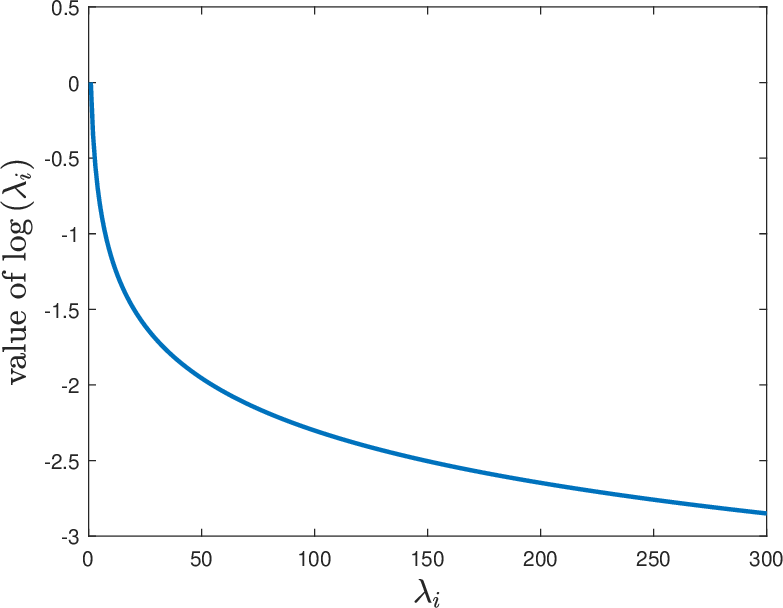}}	
		\end{center}
	\caption{Illustration of eigenvalue decay patterns in the synthetic Hessian matrices.}
	\label{fig:distr}
\end{figure*}

Figure~\ref{fig:distr} illustrates the eigenvalue distributions for the three settings. In the exponential decay case, the logarithm of eigenvalues $\log \lambda_i$ exhibits a linear decay. For the polynomial decay cases, $\log \lambda_i$ decays rapidly at first and then slows down, indicating a heavier tail in the spectrum.


\begin{figure*}[!ht]
	\begin{center}
		\centering
		\subfigure[\textsf{Convergence Rates for ``exponential decay''.}]{\includegraphics[width=55mm]{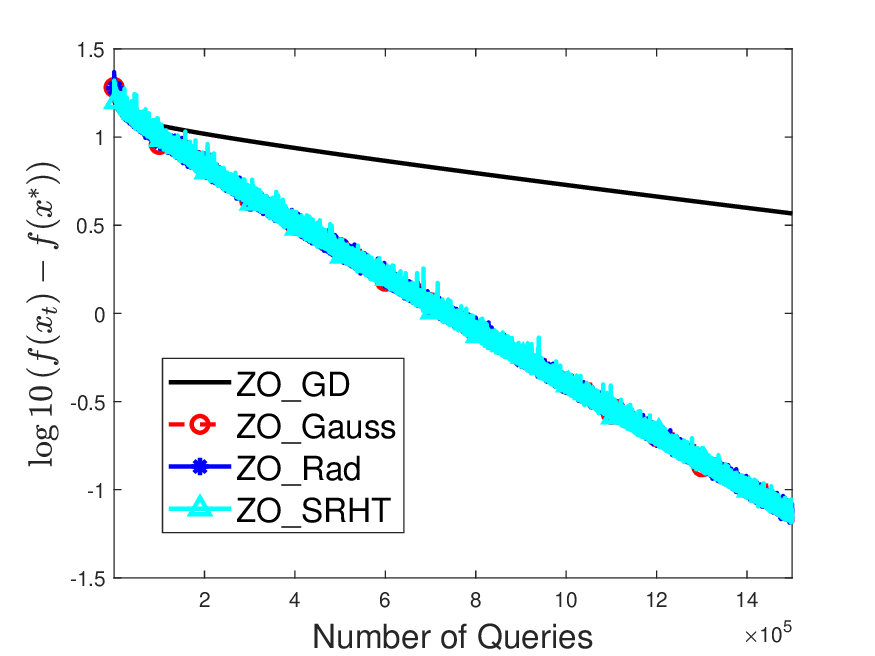}}~
		\subfigure[\textsf{Convergence Rates for ``poly decay''.}]{\includegraphics[width=55mm]{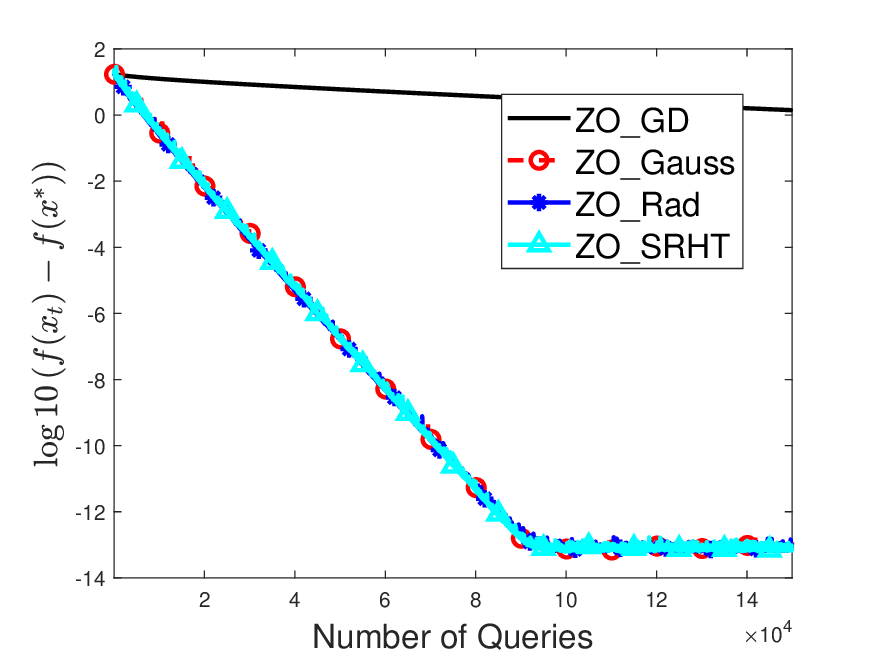}}~
		\subfigure[\textsf{Convergence Rates for ``poly decay sqrt''}]{\includegraphics[width=55mm]{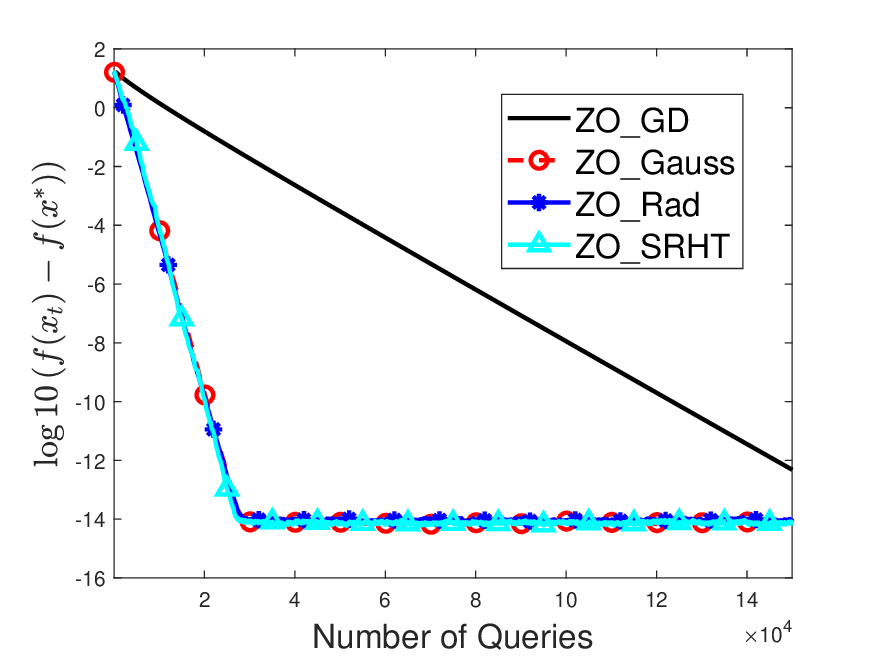}}	
		\end{center}
	\caption{Comparison of query complexity for different ZOO methods on convex quadratic functions with various Hessian eigenvalue structures. }
	\label{fig:result}
\end{figure*}


We implement Algorithm~\ref{alg:SA} with three types of sketching matrices: Gaussian, Rademacher, and SRHT, denoted as \texttt{ZO\_Gauss}, \texttt{ZO\_Rad}, and \texttt{ZO\_SRHT}, respectively. 
In all implementations, we set the sketching size $\ell = 10$.
For comparison, we also implement the FD method, which computes $d$ partial derivatives using zeroth-order oracles. 
Since FD essentially corresponds to a zeroth-order variant of gradient descent, we denote it as \texttt{ZO\_GD}.
In our experiments, we set $\alpha = 0.1$ for gradient estimation in Eq.~\eqref{eq:g}. 
Following our theoretical analysis, we choose the step size $\eta = \frac{\ell}{\tr(A) + d\lambda}$ for \texttt{ZO\_Gauss}, \texttt{ZO\_Rad}, and \texttt{ZO\_SRHT}, and set $\eta = \frac{1}{\lambda_{\max}(A) + \lambda}$ for \texttt{ZO\_GD}, in accordance with the standard theory of gradient descent~\citep{nesterov2003introductory}.

The experimental results are presented in Figure~\ref{fig:result}. 
We observe that Algorithm~\ref{alg:SA} with different random sketching matrices consistently outperforms \texttt{ZO\_GD} in all three types of Hessian matrices, each characterized by different eigenvalue distributions. 
This confirms our theoretical claim in Section~\ref{sec:convergence_convex_quad} that when $\tr(A)\ll d L$, zeroth-order methods using random sketching matrices can achieve query complexity independent of the problem dimension.
Moreover, all three types of random sketching matrices yield nearly identical query complexities.
This empirically demonstrates that using $(1/4, k, \delta/T)$-oblivious sketching matrices suffices for achieving weakly dimension-independent convergence guarantees.


\begin{table*}[]
	\centering
	\caption{Datasets summary (sparsity$=\frac{\#\text{Non-Zero Entries}}{n\times d}$)}
	\label{tb:data}
	\begin{tabular}{ccccc}
		\hline
		~~Dataset~~ &~~ $n$~~&~~$d$~~&~~$\lambda$~~&~~batch size $\ell$~~ \\ \hline
		~~a9a~~&~~$32,561$~~&~~$123$~~&~~$10^{-4}$~~&~~$10$~~  \\
		~~mushroom ~~&~~$8,124$~~&~~$112$~~&~~$10^{-4}$~~&~~$10$~~  \\ 
        ~~w8a~~ &~~ $49,749$~~&~~$300$~~&~~$10^{-4}$~~&~~$10$~~  \\
		~~gisette ~~&~~$6,000$~~&~~$5,000$~~&~~$10^{-3}$~~&~~$10$~~  \\ 
		\hline
	\end{tabular}
\end{table*}

\begin{figure*}[!ht]
	\begin{center}
		\centering
		\subfigure[\textsf{Eigenvalue distribution of the Hessian of `a9a'.}]{\includegraphics[width=55mm]{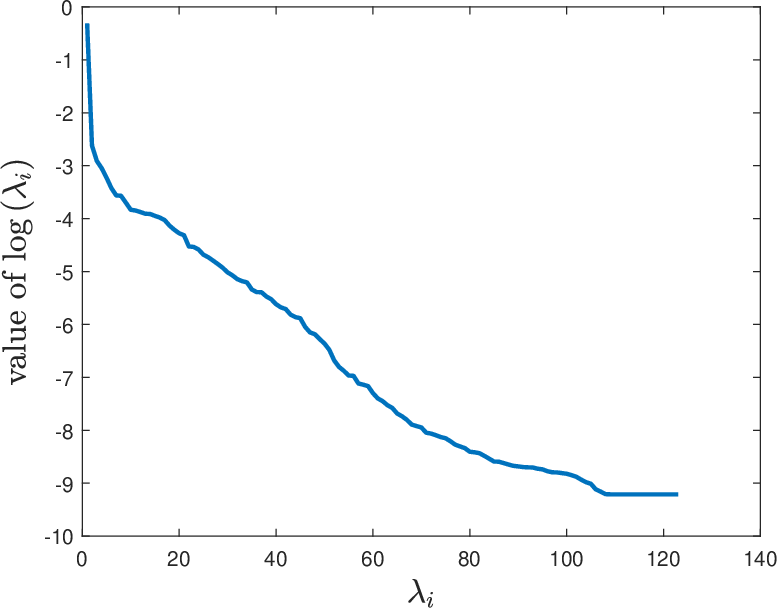}}~~~~~~~
		\subfigure[\textsf{Eigenvalue distribution of the Hessian of `mushroom' .}]{\includegraphics[width=55mm]{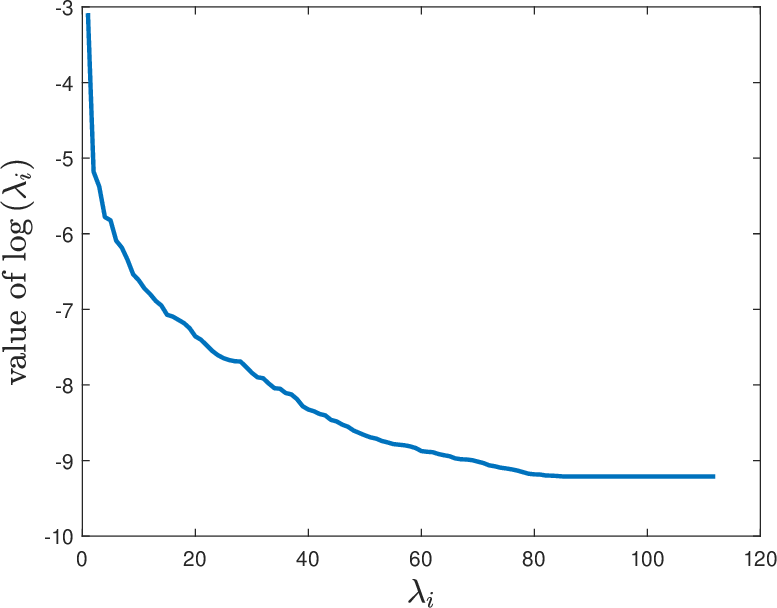}}\\
		\subfigure[\textsf{Eigenvalue distribution of the Hessian of `w8a'.}]{\includegraphics[width=55mm]{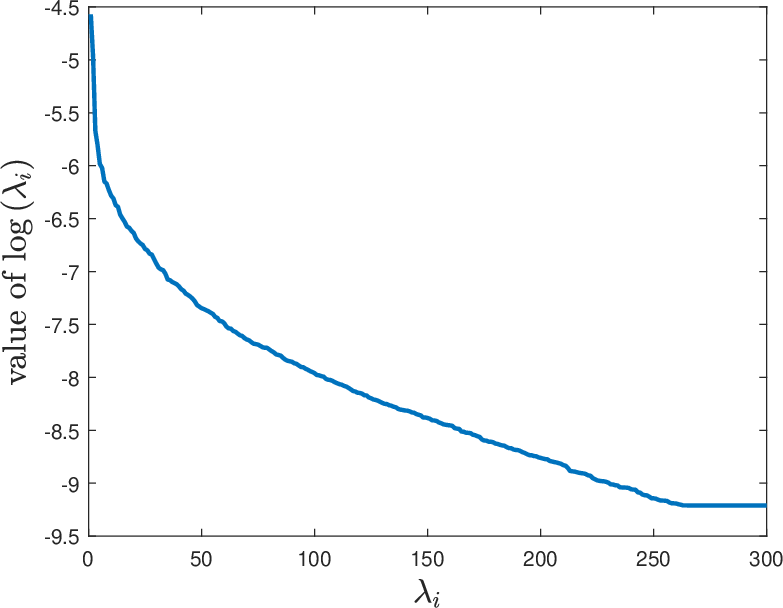}}~~~~~~~
		\subfigure[\textsf{Eigenvalue distribution of the Hessian of `gisette'.}]{\includegraphics[width=55mm]{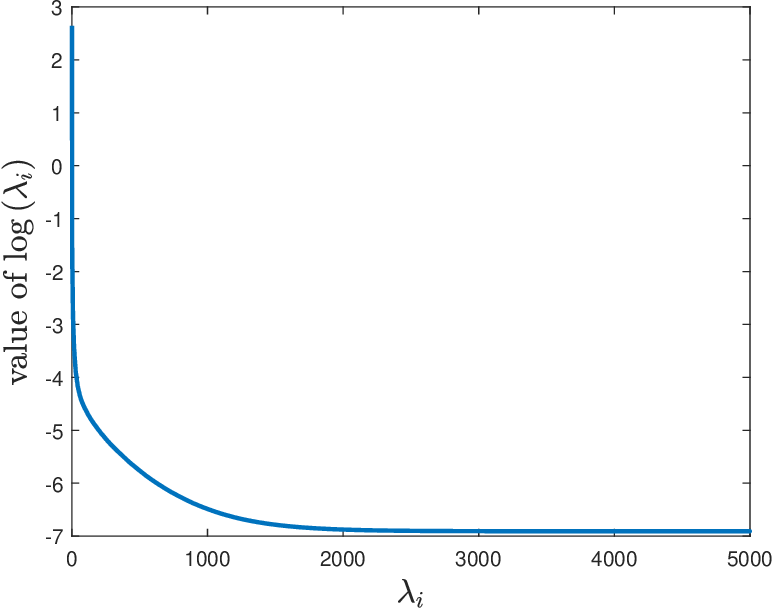}}	
	\end{center}
	\caption{Eigenvalue distributions of the Hessian matrices near the optimal solutions for four logistic regression datasets. All distributions exhibit fast decay patterns, indicating that the proposed zeroth-order algorithm is well-suited for such settings.}
	\label{fig:distr1}
\end{figure*}

\subsection{Large-Scale Logistic Regression for Real-World Classification}

Next, we evaluate our algorithm on logistic regression problems, which are non-quadratic and defined as follows:
\begin{equation}\label{eq:logi}
	\phi(\x) = \frac{1}{n}\sum_{i=1}^{n}\log\big[1+\exp(-y_i\langle \mathbf{a}_i, \x \rangle)\big] + \frac{\lambda}{2}\|\x\|^2,
\end{equation}
where $\mathbf{a}_i \in \RR^d$ denotes the $i$-th feature vector, $y_i \in \{-1, 1\}$ is the corresponding label, and $\lambda > 0$ is the regularization parameter.
We conduct experiments on four real-world datasets, that is,  ``a9a'', ``mushroom'', ``w8a'', and ``gisette''. 
The detailed description of these four datasets are listed in Table~\ref{tb:data}.
Furthermore, we set the regularization parameter $\lambda = 10^{-4}$ for ``w8a'', ``a9a'', ``mushroom'', and $\lambda = 10^{-3}$ for ``gisette''. 

To validate the applicability of our theory beyond the quadratic setting, we apply Algorithm~\ref{alg:SA} using the trace estimation strategy introduced in Section~\ref{sec:trace_estimation}. Specifically, we set the step size as $\eta_t = \frac{1}{4\tau(\x_t, S_t)}$, where $\tau(\x_t, S_t)$ is computed using Eq.~\eqref{eq:H}. For all datasets, we fix the sketching batch size to $\ell = 10$ and set $\alpha = 0.01$ for the gradient approximation in Eq.~\eqref{eq:g}.
For comparison, we also implement the \texttt{ZO\_GD} method, which computes $d$-dimensional finite difference approximations with the same $\alpha = 0.01$. The step size for \texttt{ZO\_GD} is tuned to be optimal while adhering to the classical gradient descent theory \citep{nesterov2003introductory}.

Before presenting our experimental results, we first analyze the spectral properties of the Hessian matrices near the optimal points for all four datasets. Specifically, we plot the eigenvalue distributions of the Hessians in Figure~\ref{fig:distr1}. 
We observe that all datasets exhibit a ``fast decay'' behavior in their Hessian eigenvalue spectra. Notably, the eigenvalue distributions of the ``a9a'' and ``w8a'' datasets closely resemble the ``exponential decay'' pattern. The ``gisette'' dataset shows the most rapid decay, aligning with the ``poly decay'' profile.
These observations confirm that many real-world problems naturally exhibit fast eigenvalue decay in their Hessians. Consequently, our proposed algorithm is particularly well-suited for such scenarios, where dimension-independent convergence can be effectively realized.

\begin{figure*}[!ht]
	\begin{center}
		\centering
		\subfigure[\textsf{`a9a'.}]{\includegraphics[width=55mm]{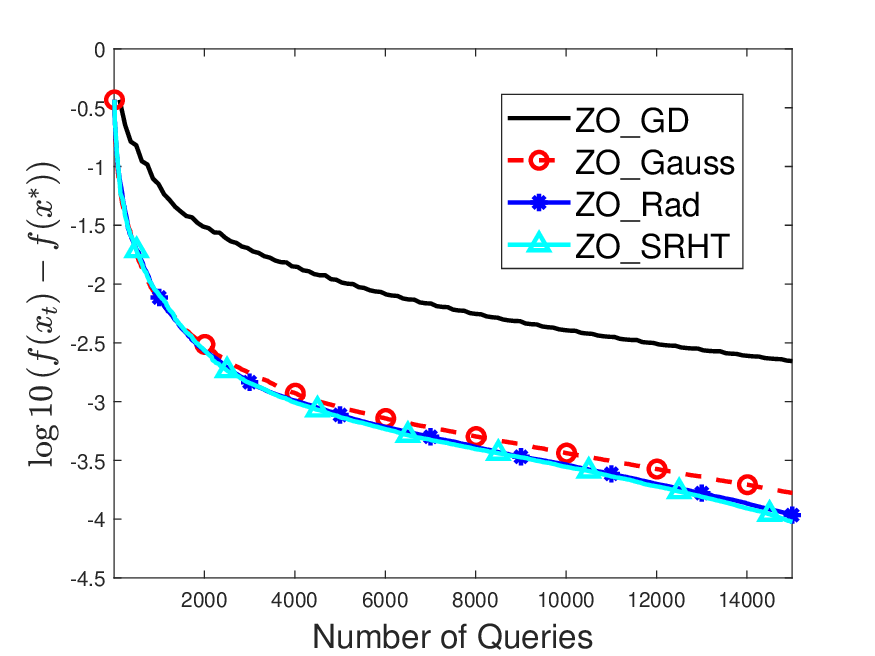}}~
		\subfigure[\textsf{`mushroom' .}]{\includegraphics[width=55mm]{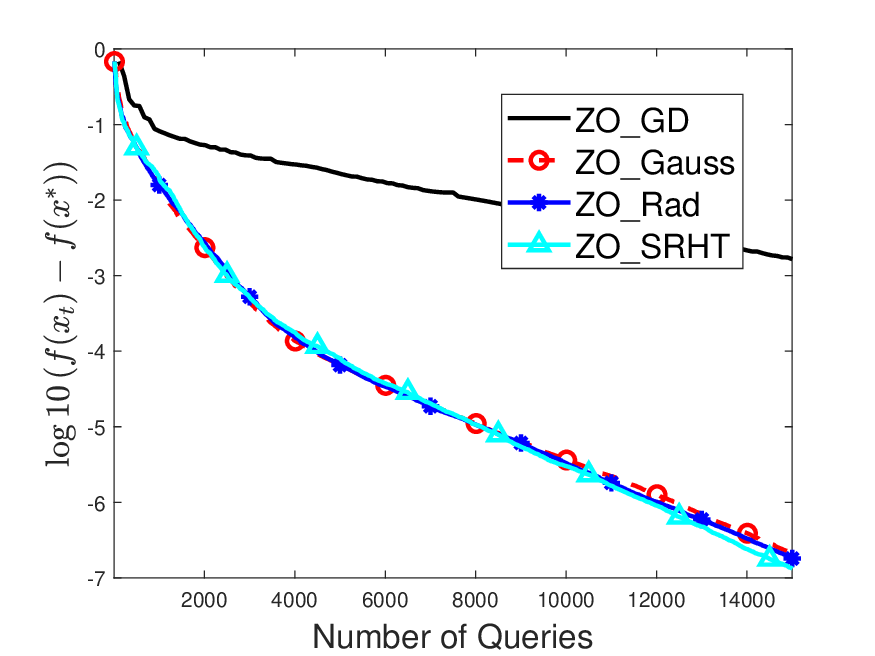}}\\
		\subfigure[\textsf{`w8a'.}]{\includegraphics[width=55mm]{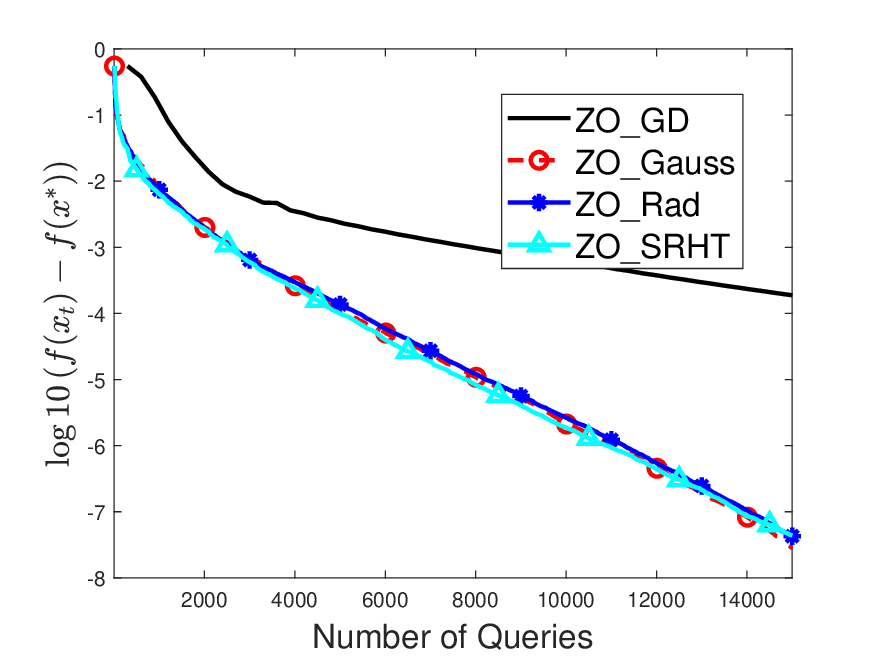}}	~
		\subfigure[\textsf{`gisette'.}]{\includegraphics[width=55mm]{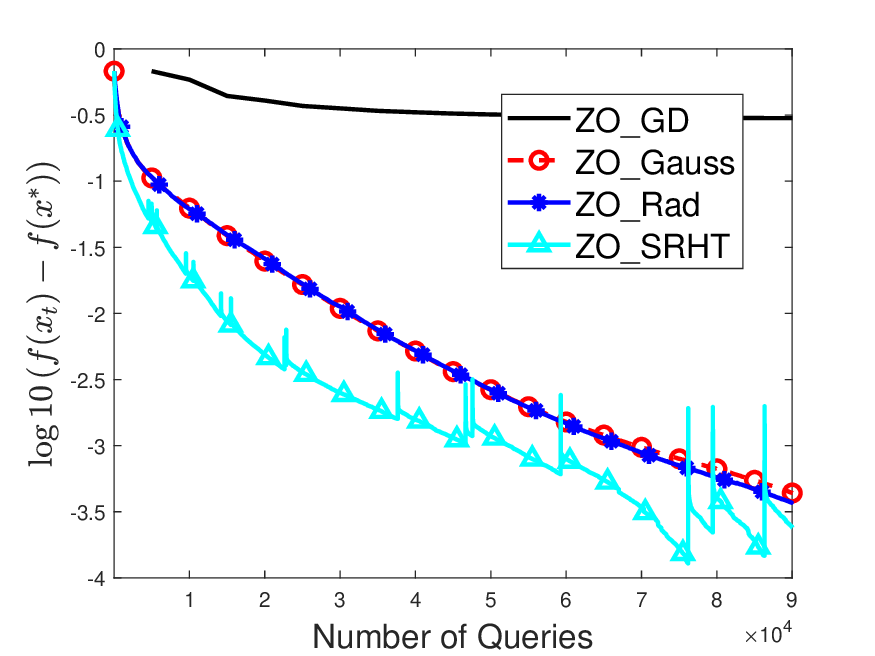}}	
	\end{center}
	\caption{Query complexity comparison of zeroth-order algorithms on real-world logistic regression datasets.}
	\label{fig:result1}
\end{figure*}

We report the experimental results in Figure~\ref{fig:result1}. 
It can be observed that our zeroth-order algorithms using random sketching matrices consistently achieve lower query complexities than \texttt{ZO\_GD}, which estimates gradients by computing $d$-partial derivatives via zeroth-order oracles.
In conjunction with the Hessian eigenvalue distributions presented in Figure~\ref{fig:distr1}, these results validate our theoretical findings in Section~\ref{sec:convegence_smooth_convex}.
Specifically, when $\tr(\nabla^2\phi(\x)) \leq dL$, our proposed zeroth-order algorithms with sketching achieve lower query complexities compared to \texttt{ZO\_GD}.
Consistent with our experiments on quadratic functions, we observe that Algorithm~\ref{alg:SA} achieves comparable performance across all three types of sketching matrices.
This reinforces the claim that $(1/4, k, \delta/T)$-oblivious sketching matrices suffice for attaining weakly dimension-independent query complexities.
Furthermore, we note that the most significant performance improvement over \texttt{ZO\_GD} occurs on the ``gisette'' dataset.
This is attributed to its Hessian exhibiting the fastest eigenvalue decay, resulting in the smallest ratio $\tr(\nabla^2\phi(\x)) / L$ among the four datasets, as shown in Figure~\ref{fig:distr1}.
This further corroborates the theoretical predictions made in Section~\ref{sec:convegence_smooth_convex}.

\section{Conclusion}

We propose a unified framework for ZOO that leverages \emph{oblivious randomized sketching} to control variance and improve query efficiency. By expressing both FD and RFD estimators through a general sketch matrix $S$, our analysis reveals that choosing $S$ from a broad class of sketching matrices—including Gaussian, Rademacher, SRHT, and sparse embeddings—enables high-probability convergence with query complexity scaling as $\tilde{\mathcal{O}}\left(\frac{\tr(A)}{\mu}\log\frac{1}{\epsilon}\right)$ under mild sketch-size conditions for convex quadratic objectives.

We further extend this framework to a Hessian-aware setting, where curvature is incorporated through sketched preconditioners without explicitly constructing Hessian blocks. The resulting algorithms retain linear per-iteration query cost and achieve improved complexity bounds of the form $\Tilde{\cO}\!\bigl(\max_t\{\rho_t^{-1}\tr(H_t^{-1}A)\}\log(1/\epsilon)\bigr)$, outperforming prior methods that require full gradient or Hessian approximations.

Overall, this work provides a general and practical foundation for designing scalable zeroth-order algorithms with theoretical guarantees. The sketch-based perspective clarifies the role of variance, dimensionality, and structure, and opens up new directions for communication-efficient, structure-aware black-box optimization in large-scale applications.


%
%
%

\bibliographystyle{informs2014}

\let\oldbibliography\thebibliography
\renewcommand{\thebibliography}[1]{%
		\oldbibliography{#1}%
		\baselineskip12pt 
		\setlength{\itemsep}{6pt}
	}
\bibliography{reference}

\begin{thebibliography}{45}
\providecommand{\natexlab}[1]{#1}
\providecommand{\url}[1]{\texttt{#1}}
\providecommand{\urlprefix}{URL }

\bibitem[{Avron \protect\BIBand{} Toledo(2011)}]{avron2011randomized}
Avron H, Toledo S (2011) Randomized algorithms for estimating the trace of an
  implicit symmetric positive semi-definite matrix. \emph{Journal of the ACM}
  58(2):1--34.

\bibitem[{Berahas et~al.(2022)Berahas, Cao, Choromanski, \protect\BIBand{}
  Scheinberg}]{berahas2022theoretical}
Berahas AS, Cao L, Choromanski K, Scheinberg K (2022) A theoretical and
  empirical comparison of gradient approximations in derivative-free
  optimization. \emph{Foundations of Computational Mathematics} 22(2):507--560.

\bibitem[{Boutsidis \protect\BIBand{} Gittens(2013)}]{boutsidis2013improved}
Boutsidis C, Gittens A (2013) Improved matrix algorithms via the subsampled
  randomized hadamard transform. \emph{SIAM Journal on Matrix Analysis and
  Applications} 34(3):1301--1340.

\bibitem[{Bubeck et~al.(2012)Bubeck, Cesa-Bianchi et~al.}]{bubeck2012regret}
Bubeck S, Cesa-Bianchi N, et~al. (2012) Regret analysis of stochastic and
  nonstochastic multi-armed bandit problems. \emph{Foundations and
  Trends{\textregistered} in Machine Learning} 5(1):1--122.

\bibitem[{Cartis \protect\BIBand{} Roberts(2023)}]{cartis2023scalable}
Cartis C, Roberts L (2023) Scalable subspace methods for derivative-free
  nonlinear least-squares optimization. \emph{Mathematical Programming}
  199(1):461--524.

\bibitem[{Choromanski et~al.(2018)Choromanski, Rowland, Sindhwani, Turner,
  \protect\BIBand{} Weller}]{choromanski2018structured}
Choromanski K, Rowland M, Sindhwani V, Turner R, Weller A (2018) Structured
  evolution with compact architectures for scalable policy optimization.
  \emph{International Conference on Machine Learning}, 970--978 (PMLR).

\bibitem[{Clarkson \protect\BIBand{} Woodruff(2017)}]{clarkson2017low}
Clarkson KL, Woodruff DP (2017) Low-rank approximation and regression in input
  sparsity time. \emph{Journal of the ACM} 63(6):1--45.

\bibitem[{Cohen(2016)}]{cohen2016simpler}
Cohen MB (2016) Simpler and tighter analysis of sparse oblivious subspace
  embeddings. \emph{Proceedings of the 27th Annual ACM-SIAM Symposium on
  Discrete Algorithms (SODA)}.

\bibitem[{Cohen et~al.(2016)Cohen, Nelson, \protect\BIBand{}
  Woodruff}]{cohen2016optimal}
Cohen MB, Nelson J, Woodruff DP (2016) Optimal approximate matrix product in
  terms of stable rank. \emph{International Colloquium on Automata, Languages,
  and Programming} (Schloss Dagstuhl-Leibniz-Zentrum fur Informatik GmbH,
  Dagstuhl Publishing).

\bibitem[{Conn et~al.(1997)Conn, Scheinberg, \protect\BIBand{}
  Toint}]{conn1997convergence}
Conn AR, Scheinberg K, Toint PL (1997) On the convergence of derivative-free
  methods for unconstrained optimization. \emph{Approximation Theory and
  Tptimization: tributes to MJD Powell} 1:83--108.

\bibitem[{Conn et~al.(2009)Conn, Scheinberg, \protect\BIBand{}
  Vicente}]{conn2009introduction}
Conn AR, Scheinberg K, Vicente LN (2009) \emph{Introduction to derivative-free
  optimization} (SIAM).

\bibitem[{Cortinovis \protect\BIBand{}
  Kressner(2021)}]{cortinovis2021randomized}
Cortinovis A, Kressner D (2021) On randomized trace estimates for indefinite
  matrices with an application to determinants. \emph{Foundations of
  Computational Mathematics} 1--29.

\bibitem[{Deng(2012)}]{deng2012mnist}
Deng L (2012) The mnist database of handwritten digit images for machine
  learning research. \emph{IEEE Signal Processing Magazine} 29(6):141--142.

\bibitem[{Ghadimi \protect\BIBand{} Lan(2013)}]{ghadimi2013stochastic}
Ghadimi S, Lan G (2013) Stochastic first-and zeroth-order methods for nonconvex
  stochastic programming. \emph{SIAM Journal on Optimization} 23(4):2341--2368.

\bibitem[{Gower et~al.(2019)Gower, Kovalev, Lieder, \protect\BIBand{}
  Richt{\'a}rik}]{gower2019rsn}
Gower R, Kovalev D, Lieder F, Richt{\'a}rik P (2019) \text{RSN}: randomized
  subspace newton. \emph{Advances in Neural Information Processing Systems} 32.

\bibitem[{Hansen(2016)}]{hansen2016cma}
Hansen N (2016) The \text{CMA} evolution strategy: A tutorial. \emph{arXiv
  preprint arXiv:1604.00772} .

\bibitem[{Ilyas et~al.(2018)Ilyas, Engstrom, Athalye, \protect\BIBand{}
  Lin}]{ilyas2018black}
Ilyas A, Engstrom L, Athalye A, Lin J (2018) Black-box adversarial attacks with
  limited queries and information. \emph{International conference on machine
  learning}, 2137--2146 (PMLR).

\bibitem[{James et~al.(2013)James, Witten, Hastie, Tibshirani
  et~al.}]{james2013introduction}
James G, Witten D, Hastie T, Tibshirani R, et~al. (2013) \emph{An introduction
  to statistical learning}, volume 112 (Springer).

\bibitem[{Kiefer \protect\BIBand{} Wolfowitz(1952)}]{kiefer1952stochastic}
Kiefer J, Wolfowitz J (1952) Stochastic estimation of the maximum of a
  regression function. \emph{The Annals of Mathematical Statistics} 462--466.

\bibitem[{Lam \protect\BIBand{} Zhang(2024)}]{lam2024distributionally}
Lam H, Zhang J (2024) Distributionally constrained black-box stochastic
  gradient estimation and optimization. \emph{Operations Research}
  75(5):2680--2694.

\bibitem[{Ma \protect\BIBand{} Huang(2025)}]{ma2025revisiting}
Ma S, Huang H (2025) Revisiting zeroth-order optimization: Minimum-variance
  two-point estimators and directionally aligned perturbations. \emph{The
  Thirteenth International Conference on Learning Representations}.

\bibitem[{Malladi et~al.(2023)Malladi, Gao, Nichani, Damian, Lee, Chen,
  \protect\BIBand{} Arora}]{malladi2023fine}
Malladi S, Gao T, Nichani E, Damian A, Lee JD, Chen D, Arora S (2023)
  Fine-tuning language models with just forward passes. \emph{Advances in
  Neural Information Processing Systems} 36:53038--53075.

\bibitem[{Meng \protect\BIBand{} Mahoney(2013)}]{meng2013low}
Meng X, Mahoney MW (2013) Low-distortion subspace embeddings in input-sparsity
  time and applications to robust linear regression. \emph{Proceedings of the
  forty-fifth annual ACM symposium on Theory of computing}, 91--100.

\bibitem[{Meyer et~al.(2021)Meyer, Musco, Musco, \protect\BIBand{}
  Woodruff}]{meyer2021hutch}
Meyer RA, Musco C, Musco C, Woodruff DP (2021) Hutch++: Optimal stochastic
  trace estimation. \emph{Symposium on Simplicity in Algorithms (SOSA)},
  142--155 (SIAM).

\bibitem[{Nelson \protect\BIBand{} Nguy{\^e}n(2013)}]{nelson2013osnap}
Nelson J, Nguy{\^e}n HL (2013) Osnap: Faster numerical linear algebra
  algorithms via sparser subspace embeddings. \emph{2013 ieee 54th annual
  symposium on foundations of computer science}, 117--126 (IEEE).

\bibitem[{Nesterov(2003)}]{nesterov2003introductory}
Nesterov Y (2003) \emph{Introductory lectures on convex optimization: A basic
  course}, volume~87 (Springer Science \& Business Media).

\bibitem[{Nesterov(2012)}]{nesterov2012efficiency}
Nesterov Y (2012) Efficiency of coordinate descent methods on huge-scale
  optimization problems. \emph{SIAM Journal on Optimization} 22(2):341--362.

\bibitem[{Nesterov \protect\BIBand{} Spokoiny(2017)}]{nesterov2017random}
Nesterov Y, Spokoiny V (2017) Random gradient-free minimization of convex
  functions. \emph{Foundations of Computational Mathematics} 17:527--566.

\bibitem[{Powell(1972)}]{powell1972unconstrained}
Powell MJ (1972) Unconstrained minimization algorithms without computation of
  derivatives. Technical report.

\bibitem[{Roosta-Khorasani \protect\BIBand{} Ascher(2015)}]{roosta2015improved}
Roosta-Khorasani F, Ascher U (2015) Improved bounds on sample size for implicit
  matrix trace estimators. \emph{Foundations of Computational Mathematics}
  15(5):1187--1212.

\bibitem[{Sagun et~al.(2016)Sagun, Bottou, \protect\BIBand{}
  LeCun}]{sagun2016eigenvalues}
Sagun L, Bottou L, LeCun Y (2016) Eigenvalues of the hessian in deep learning:
  Singularity and beyond. \emph{arXiv preprint arXiv:1611.07476} .

\bibitem[{Scheinberg(2022)}]{scheinberg2022finite}
Scheinberg K (2022) Finite difference gradient approximation: To randomize or
  not? \emph{INFORMS Journal on Computing} 34(5):2384--2388.

\bibitem[{Spall(1992)}]{spall1992multivariate}
Spall JC (1992) Multivariate stochastic approximation using a simultaneous
  perturbation gradient approximation. \emph{IEEE Transactions on Automatic
  Control} 37(3):332--341.

\bibitem[{Wang(2025)}]{wang2025adaptivity}
Wang Y (2025) On adaptivity in nonstationary stochastic optimization with
  bandit feedback. \emph{Operations Research} 73(2):819--828.

\bibitem[{Wang et~al.(2024)Wang, Ye, Dai, \protect\BIBand{} Tsang}]{wangcan}
Wang Y, Ye H, Dai G, Tsang I (2024) Can {G}aussian sketching converge faster on
  a preconditioned landscape? \emph{Proceedings of the 41st International
  Conference on Machine Learning}, volume 235, 52064--52082 (PMLR).

\bibitem[{Wierstra et~al.(2014)Wierstra, Schaul, Glasmachers, Sun, Peters,
  \protect\BIBand{} Schmidhuber}]{wierstra2014natural}
Wierstra D, Schaul T, Glasmachers T, Sun Y, Peters J, Schmidhuber J (2014)
  Natural evolution strategies. \emph{The Journal of Machine Learning Research}
  15(1):949--980.

\bibitem[{Wild et~al.(2008)Wild, Regis, \protect\BIBand{}
  Shoemaker}]{wild2008orbit}
Wild SM, Regis RG, Shoemaker CA (2008) Orbit: Optimization by radial basis
  function interpolation in trust-regions. \emph{SIAM Journal on Scientific
  Computing} 30(6):3197--3219.

\bibitem[{Woodruff et~al.(2014)}]{woodruff2014sketching}
Woodruff DP, et~al. (2014) Sketching as a tool for numerical linear algebra.
  \emph{Foundations and Trends{\textregistered} in Theoretical Computer
  Science} 10(1--2):1--157.

\bibitem[{Yao et~al.(2020)Yao, Gholami, Keutzer, \protect\BIBand{}
  Mahoney}]{yao2020pyhessian}
Yao Z, Gholami A, Keutzer K, Mahoney MW (2020) Pyhessian: Neural networks
  through the lens of the hessian. \emph{2020 IEEE International Conference on
  Big Data}, 581--590 (IEEE).

\bibitem[{Ye(2023)}]{ye2023mirror}
Ye H (2023) Mirror natural evolution strategies. \emph{arXiv preprint
  arXiv:2308.00469} .

\bibitem[{Ye et~al.(2025)Ye, Huang, Fang, Li, \protect\BIBand{}
  Zhang}]{ye2018hessian}
Ye H, Huang Z, Fang C, Li CJ, Zhang T (2025) Hessian-aware zeroth-order
  optimization. \emph{IEEE Transactions on Pattern Analysis and Machine
  Intelligence (early access)} .

\bibitem[{Ye et~al.(2021)Ye, Luo, \protect\BIBand{} Zhang}]{ye2021approximate}
Ye H, Luo L, Zhang Z (2021) Approximate newton methods. \emph{The Journal of
  Machine Learning Research} 22(1):3067--3107.

\bibitem[{Yue et~al.(2023{\natexlab{a}})Yue, Fang, \protect\BIBand{}
  Lin}]{yue2023lower}
Yue P, Fang C, Lin Z (2023{\natexlab{a}}) On the lower bound of minimizing
  polyak-{\l}ojasiewicz functions. \emph{The Thirty Sixth Annual Conference on
  Learning Theory}, 2948--2968 (PMLR).

\bibitem[{Yue et~al.(2023{\natexlab{b}})Yue, Yang, Fang, \protect\BIBand{}
  Lin}]{yue2023zeroth}
Yue P, Yang L, Fang C, Lin Z (2023{\natexlab{b}}) Zeroth-order optimization
  with weak dimension dependency. \emph{The Thirty Sixth Annual Conference on
  Learning Theory}, 4429--4472 (PMLR).

\bibitem[{Zhang et~al.(2024)Zhang, Thekumparampil, Oh, \protect\BIBand{}
  He}]{zhang2024dpzero}
Zhang L, Thekumparampil KK, Oh S, He N (2024) Dpzero: dimension-independent and
  differentially private zeroth-order optimization (International Conference on
  Machine Learning (ICML 2024)).

\end{thebibliography}

\newpage
\begin{APPENDICES}

\section{Proof of Lemmas}

\subsection{Proof of Lemma \ref{lem:g_prop}}
By Assumption~\ref{ass:phi}, we have
$\phi(\x + \alpha \si) = \phi(\x) + \alpha \dotprod{\nabla \phi(\x), \si } +  \frac{\alpha^2}{2} \norm{\si}_A^2$, and
$\phi(\x - \alpha \si) = \phi(\x) - \alpha \dotprod{\nabla \phi(\x), \si } +  \frac{\alpha^2}{2} \norm{\si}_A^2$.
Thus, we have
\begin{equation}\label{eq:phi_diff}
\frac{\phi(\x + \alpha \si) - \phi(\x - \alpha \si)}{2\alpha} = \dotprod{\nabla \phi(\x), \si}.
\end{equation}
Thus, we can obtain that
\begin{align*}
	&\g(\x) 
	=\sum_{i=1}^\ell \frac{f(\x+\alpha \si) - f(\x - \alpha \si)}{2\alpha} \si\\
	&=
	\sum_{i=1}^{\ell} \frac{\phi(\x + \alpha \si) - \phi(\x - \alpha \si) }{2\alpha} \si+\sum_{i=1}^{\ell}\frac{\zeta(\x + \alpha \si) - \zeta(\x - \alpha \si)}{2\alpha} \si\\
	&\stackrel{\eqref{eq:phi_diff}}{=}
	\sum_{i=1}^{\ell}\dotprod{\nabla\phi(\x), \si} \si+  \sum_{i=1}^{\ell} \frac{\zeta(\x + \alpha \si) - \zeta(\x - \alpha \si)}{2\alpha} \si\\
	&=
	SS^\top \nabla \phi(\x) + S \bv.
\end{align*}
By Assumption~\ref{ass:zeta}, we can obtain that $\norm{\bv}^2 \le \frac{\ell\sigma^2}{\alpha^2}$ which concludes the proof.

\subsection{Proof of Lemma \ref{lemma:bound_phi}}

First, using the property of the quadratic function, we have
\begin{align*}
	&\phi(\x_{t+1}) \\ 
	\stackrel{\eqref{eq:update}\eqref{eq:phi}}{=}&
	\phi(\x_t) - \eta \dotprod{\nabla \phi(\x_t), \g(\x_t)} + \frac{\eta^2}{2} \norm{\g(\x_t)}_A^2\\
	\stackrel{\eqref{eq:g_prop}}{=}&
	\phi(\x_t) - \eta \dotprod{\nabla \phi(\x_t), S_tS_t^\top \nabla \phi(\x_t)} - \eta \dotprod{\nabla \phi(\x_t), S_t\bv_t} + \frac{\eta^2}{2} \norm{S_tS_t^\top \nabla \phi(\x_t) + S_t\bv_t}_A^2\\
	\le&
	\phi(\x_t) - \eta \norm{S_t^\top \nabla \phi(\x_t)}^2 + \frac{\eta}{2} \norm{S_t^\top\nabla \phi(\x_t)}^2 + \frac{\eta}{2} \norm{\bv_t}^2 + \eta^2 \norm{S_tS_t^\top \nabla \phi(\x_t)}_A^2 + \eta^2\norm{S_t \bv_t}_A^2\\
	=&
	\phi(\x_t) - \frac{\eta}{2} \norm{S_t^\top \nabla \phi(\x_t)}^2 + \frac{\eta}{2} \norm{\bv_t}^2 + \eta^2 \norm{S_tS_t^\top \nabla \phi(\x_t)}_A^2 + \eta^2\norm{S_t \bv_t}_A^2,
\end{align*}
where the first inequality is because of the Cauchy's inequality and $2ab \le a^2 + b^2$.

Furthermore, 
\begin{align*}
\norm{S_tS_t^\top \nabla \phi(\x_t)}_A^2 
=
\nabla^\top \phi(\x_t) S_tS_t^\top A S_t S_t^\top \nabla\phi(\x_t)
\le
\norm{S_t^\top A S_t}_2 \cdot \norm{S_t^\top \nabla\phi(\x_t)}^2.
\end{align*}
Similarly, we can obtain that
\begin{align*}
\norm{S_t \bv_t}_A^2 
= 
\bv_t^\top S_t^\top A S_t \bv_t 
\le
\norm{S_t^\top A S_t}_2 \cdot \norm{\bv_t}^2.
\end{align*}

Therefore, we can obtain that
\begin{align*}
	\phi(\x_{t+1}) 
	\leq
	\phi(\x_t) - \frac{\eta}{2} \norm{S_t^\top \nabla\phi(\x_t)}^2 + \eta^2 \norm{S_t^\top A S_t}_2 \cdot \norm{S_t^\top \nabla \phi(\x_t)}^2 + \frac{\eta}{2} \norm{\bv_t}^2 + \eta^2 \norm{S_t^\top A S_t}_2 \cdot \norm{\bv_t}^2. 
\end{align*}

\subsection{Proof of Lemma \ref{lem:ske}}

By the fact that $A$ is positive definite, we have $A = A^{1/2} A^{1/2}$. 
	Combining the fact that  $\norm{B^\top B}_2 = \norm{B B^\top}_2$ holds for all matrices, we have
\begin{equation}\label{eq:sas}
	\norm{S_t^\top A S_t}_2 
	= 
	\norm{S_t^\top A^{1/2} A^{1/2} S_t}_2
	=
	\norm{A^{1/2} S_t S_t^\top A^{1/2}}_2.
\end{equation}
Using the condition that $S_t$ is a $(\frac{1}{4}, k, \delta)$-oblivious sketching matrix for the matrix product, we have the following property holding with a probability at least $1-\delta$,
\begin{align*}
	\norm{A^{1/2} S_t S_t^\top A^{1/2} - A^{1/2} A^{1/2}}_2 
	\le 
	\frac{1}{4}\left(\norm{A^{1/2}}_2^2 + \frac{\norm{A^{1/2}}_F^2}{k}\right)
	\stackrel{\eqref{eq:tr_F}}{=}
	\frac{\norm{A}_2}{4} + \frac{\tr(A)}{4k}.
\end{align*}
Using the triangle inequality, we have
\begin{align*}
	\norm{A^{1/2} S_t S_t^\top A^{1/2} - A^{1/2} A^{1/2}}_2 
	\ge
	\norm{A^{1/2} S_t S_t^\top A^{1/2}}_2 - \norm{A^{1/2} A^{1/2}}_2
	\stackrel{\eqref{eq:sas}}{=} 
	\norm{S_t^\top A S_t}_2 - \norm{A}_2.
\end{align*}
Combining above two inequalities, we can obtain that
\begin{align*}
\norm{S_t^\top A S_t}_2 
\le
\frac{5\norm{A}_2}{4} + \frac{\tr(A)}{4k}.
\end{align*}

Similarly, by the condition that $S_t$ is a $(\frac{1}{4}, k, \delta)$-oblivious sketching matrix for the matrix product, we have
\begin{align*}
\norm{\nabla^\top \phi(\x_t)S_t S_t^\top \nabla \phi(x_t) - \nabla^\top \phi(\x_t) \nabla \phi(\x_t) }_2 \le \frac{1}{4} \left( \norm{\nabla \phi(\x_t)}^2 + \frac{\norm{\nabla \phi(\x_t)}^2}{k} \right) 
\le \frac{1}{2} \norm{\nabla \phi(\x_t)}^2.
\end{align*}
Using the triangle inequality, we have
\begin{align*}
\norm{\nabla^\top \phi(\x_t)S_t S_t^\top \nabla \phi(\x_t) - \nabla^\top \phi(\x_t) \nabla \phi(\x_t) }_2 
\ge& 
\norm{\nabla \phi(\x_t)}^2 - \norm{\nabla^\top \phi(\x_t)S_t S_t^\top \nabla \phi(\x_t)}\\
=&
\norm{\nabla \phi(\x_t)}^2 - \norm{S_t^\top \nabla \phi(\x_t)}^2.
\end{align*}
Similarly, we can obtain that
\begin{align*}
\norm{\nabla^\top \phi(\x_t)S_t S_t^\top \nabla \phi(\x_t) - \nabla^\top \phi(\x_t) \nabla \phi(\x_t) }_2 
\ge 
\norm{S_t^\top \nabla \phi(\x_t)}^2 - \norm{\nabla \phi(\x_t)}^2.
\end{align*}
Combining with above two inequalities, we have
\begin{align*}
\frac{1}{2} \norm{\nabla \phi(\x_t)}^2 \leq 	\norm{S_t^\top \nabla \phi(\x_t)}^2 \leq \frac{3}{2}\norm{\nabla \phi(\x_t)}^2.
\end{align*}

\subsection{Proof of Lemma \ref{lem:tg_prop}}

By Assumption~\ref{ass:phi}, we have
\begin{align*}
	\phi(\x_t + \alpha \sti) =& \phi(\x_t) + \alpha \dotprod{\nabla \phi(\x_t), H_t^{-1/2}\sti } +  \frac{\alpha^2}{2} \norm{H_t^{-1/2}\sti}_A^2,\\
	\phi(\x_t - \alpha \sti) =& \phi(\x_t) - \alpha \dotprod{\nabla \phi(\x_t), H_t^{-1/2}\sti } +  \frac{\alpha^2}{2} \norm{H_t^{-1/2}\sti}_A^2.
\end{align*}
Thus, we have
\begin{equation}\label{eq:phi_diff1}
	\frac{\phi(\x_t + \alpha H_t^{-1/2}\sti) - \phi(\x_t - \alpha H_t^{-1/2}\sti)}{2\alpha} = \dotprod{\nabla \phi(\x_t), H_t^{-1/2}\sti}.
\end{equation}
Thus, we can obtain that
\begin{align*}
	\tg(\x_t) 
	=& 
	\sum_{i=1}^\ell \frac{f(\x_t+\alpha H_t^{-1/2}\sti) - f(\x_t - \alpha H_t^{-1/2}\sti)}{2\alpha} H_t^{-1/2}\sti\\
	=&
	\sum_{i=1}^{\ell} \frac{\phi(\x_t + \alpha H_t^{-1/2}\sti) - \phi(\x_t - \alpha H_t^{-1/2}\sti) + \zeta(\x_t + \alpha H_t^{-1/2}\sti) - \zeta(\x_t - \alpha H_t^{-1/2}\sti)}{2\alpha} H_t^{-1/2}\sti \\
	\stackrel{\eqref{eq:phi_diff1}}{=}&
	\sum_{i=1}^{\ell}\dotprod{\nabla\phi(\x_t), H_t^{-1/2}\sti} H_t^{-1/2}\sti +  \sum_{i=1}^{\ell} \frac{\zeta(\x_t + \alpha H_t^{-1/2}\sti) - \zeta(\x_t - \alpha H_t^{-1/2} \sti)}{2\alpha} H_t^{-1/2}\sti\\
	=&
	H_t^{-1/2}S_tS_t^\top H_t^{-1/2}\nabla \phi(\x_t) + H_t^{-1/2}S_t \tv_t.
\end{align*}
By Assumption~\ref{ass:zeta}, we can obtain that $	\norm{\tv_t}^2 \leq \frac{\ell\sigma^2}{\alpha^2}$, which conclude the proof.

\subsection{Proof of Lemma \ref{lem:g_prop1}}

By the Taylor's expansion and Assumption~\ref{ass:L}-\ref{ass:mL}, we have
\begin{equation}\label{eq:mL}
	\left| \phi(\y) - \phi(\x) - \dotprod{\nabla \phi(\x), \y-\x} - \frac{1}{2} \norm{\y-\x}_{\nabla^2 \phi(\x)}^2  \right| \le \frac{\mL}{6}\norm{\y - \x}^3.
\end{equation}
We have
\begin{align*}
	g(\x_t) 
	=& 
	\sum_{i=1}^\ell \frac{f(\x_t+\alpha \sti) - f(\x_t - \alpha\sti)}{2\alpha}\sti\\
	=&
	\sum_{i=1}^\ell \frac{\phi(\x_t+\alpha \sti) - \phi(\x_t-\alpha \sti) + \zeta(\x_t+\alpha \sti) - \zeta(\x_t-\alpha \sti)}{2\alpha}\\
	=&
	S_tS_t^\top \nabla \phi(\x_t) + \sum_{i=1}^{\ell} \frac{\phi(\x_t+\alpha \sti) - \phi(\x_t-\alpha \sti) - 2 \alpha \dotprod{\nabla \phi(\x_t), \sti}}{2 \alpha}\sti \\
	 + & \frac{\zeta(\x_t+\alpha \sti) - \zeta(\x_t-\alpha \sti)}{2\alpha} \sti\\
	=& 
	S_tS_t^\top \nabla \phi(\x_t) + S_t \bw_t + S_t \bv_t.
\end{align*}

Letting us denote that $ D(\alpha) = \phi(\x_t + \alpha \sti)- \phi(\x_t) - \dotprod{\nabla \phi(x)_t,\alpha \sti} - \frac{\alpha^2}{2} \norm{\sti}_{\nabla^2 \phi(\x_t)}^2  $, then we have
\begin{align*}
|\bw_t^{(i)}| 
=&\left|\frac{\phi(\x_t+\alpha \sti) - \phi(\x_t-\alpha \sti) - 2 \alpha \dotprod{\nabla \phi(\x_t), \sti}}{2 \alpha}\right|
=
\left|\frac{D(\alpha) - D(-\alpha)}{2\alpha}\right| \\
\leq&
\frac{|D(\alpha)|+|D(-\alpha)|}{2\alpha}
\stackrel{\eqref{eq:mL}}{\le} \frac{\mL\alpha^2\norm{\sti}^3}{6}.
\end{align*}
Thus, 
\begin{align*}
	\norm{\bw_t}^2 = \sum_{i=1}^\ell \left(\bw_t^{(i)}\right)^2 \le \frac{\mL^2\alpha^4}{36} \sum_{i=1}^{\ell} \norm{\sti}^6 
	\leq
	\frac{\mL^2\alpha^4}{36} \norm{S_t}_F^6.
\end{align*}

\subsection{Proof of Lemma \ref{lemma:th_part1}}

	By the update rule and Assumption~\ref{ass:L}, we have
	\begin{align*}
		&\phi(\x_{t+1}) \\
		\le& 
		\phi(\x_t) - \eta_t \dotprod{\nabla \phi(\x_t), g(\x_t)} + \frac{L\eta_t^2}{2} \norm{g(\x_t)}_2^2\\
		\stackrel{\eqref{eq:g_prop1}}{=}&
		\phi(\x_t) - \eta_t \dotprod{\nabla \phi(\x_t), S_tS_t^\top \nabla \phi(\x_t) +S_t \bw_t + S_t \bv_t} + \frac{L\eta_t^2}{2} \norm{S_tS_t^\top \nabla \phi(\x_t) + S_t \bw_t + S_t \bv_t}_2^2\\
		\le&
		\phi(\x_t) - \frac{\eta_t}{2} \norm{S_t^\top\nabla \phi(\x_t)}_2^2 + \eta_t \norm{\bw_t}_2^2 + \eta_t \norm{\bv_t}_2^2  + L\eta_t^2 \norm{S_tS_t^\top}_2\left( \norm{S_t^\top \nabla \phi(\x_t)}_2^2 + 2\norm{\bw_t}_2^2 + 2\norm{\bv_t}_2^2\right)
	\end{align*}
	By Lemma~\ref{lem:ske} with $A$ being the identity matrix, we have
	\begin{align*}
		\norm{S_tS_t^\top}_2 \stackrel{\eqref{eq:sas1}}{\leq} \frac{5}{4} + \frac{d}{4k},
	\end{align*}
    with a probability at least $1-\delta/T$.
	Thus, it holds
	\begin{align*}
		\phi(\x_{t+1}) -\phi(\x^*)
		\le &
		\phi(\x_t) - \phi(\x^*) - \eta_t\left(\frac{1}{2} - L\left(\frac{5}{4} + \frac{d}{4k}\right) \eta_t \right) \cdot \norm{S_t^\top \nabla \phi(\x_t)}^2\\
		+ &\eta_t\left(1 + 2L\left(\frac{5}{4} + \frac{d}{4k}\right) \eta_t\right)\cdot \left(\norm{\bw_t}^2 + \norm{\bv_t}^2\right)\\
		=&
		\phi(\x_t) - \phi(\x^*) - \frac{1}{4L}\left(5 + \frac{d}{k}\right)^{-1}\cdot \norm{S_t^\top \nabla \phi(\x_t)}^2
		+ \frac{3}{2L} \left(5 + \frac{d}{k}\right)^{-1} \cdot \left(\norm{\bw_t}^2 + \norm{\bv_t}^2\right)\\
		\stackrel{\eqref{eq:sas1}}{\le}&
		\phi(\x_t) - \phi(\x^*) - \frac{1}{8L}\left(5 + \frac{d}{k}\right)^{-1}\cdot \norm{\nabla \phi(\x_t)}^2
		+ \frac{3}{2L}\left(5 + \frac{d}{k}\right)^{-1} \cdot \left(\norm{\bw_t}^2 + \norm{\bv_t}^2\right)\\
		\leq&
		\left(1 - \frac{\mu}{4L} \left(5 + \frac{d}{k}\right)^{-1}\right) \cdot \Big( \phi(\x_t) - \phi(\x^*) \Big) +  \frac{3}{2L}\left(5 + \frac{d}{k}\right)^{-1} \cdot \left(\frac{\mL^2\alpha^4}{36} \norm{S_t}_F^6+\frac{\ell\sigma^2}{\alpha^2}\right),
	\end{align*}
with a probability at least $1-\delta/T$.

\subsection{Proof of Lemma \ref{lemma:th_part2}}

	By the $\mL$-Hessian Lipschitz continuity, we can obtain that
\begin{align*}
	\phi(\x_{t+1})
	\stackrel{\eqref{eq:mL}}{\le}&
	\phi(\x_t) - \eta_t\dotprod{\nabla \phi(\x_t), g(\x_t)} + \frac{\eta_t^2}{2} \norm{g(\x_t)}_{\nabla^2\phi(\x_t)}^2  + \frac{\mL\eta_t^3}{6} \norm{g(\x_t)}^3\\
	\stackrel{\eqref{eq:g_prop1}}{=}&
	\phi(\x_t)- \eta_t \dotprod{\nabla \phi(\x_t), S_tS_t^\top \nabla \phi(\x_t) +S_t \bw_t + S_t \bv_t} + \frac{\eta_t^2}{2} \norm{S_tS_t^\top \nabla \phi(\x_t) + S_t \bw_t + S_t \bv_t}_{\nabla^2\phi(\x_t)}^2\\
	+&\frac{\mL\eta_t^3}{6} \norm{ S_tS_t^\top \nabla \phi(\x_t) + S_t \bw_t + S_t \bv_t }^3\\
	\le&
	\phi(\x_t) - \frac{\eta_t}{2} \norm{S_t^\top \nabla \phi(\x_t)}^2 + \eta_t \left( \norm{\bw_t}^2 + \norm{\bv_t}^2 \right) \\
	+& \eta_t^2 \cdot \norm{S_t^\top \nabla^2\phi(\x_t) S_t}_2 \cdot \left( \norm{S_t^\top \nabla\phi(\x_t)}^2 + 2\norm{\bw_t}^2 + 2\norm{\bv_t}^2\right)\\
	+& \frac{2\mL\eta_t^3}{3} \norm{S_t}_2^3 \left( \norm{S_t^\top \nabla \phi(\x_t)}^3 + 4\norm{\bw_t}^3 + 4\norm{\bv_t}^3 \right).
\end{align*}

By Lemma~\ref{lem:ske} and replacing $A$ with $\nabla^2 \phi(\x_t)$ and the identity matrix respectively, we can obtain that
\begin{align*}
\norm{S_t^\top \nabla^2\phi(\x_t) S_t}_2 
\le
\frac{5 \norm{\nabla^2 \phi(\x_t)}_2}{4} + \frac{\tr\Big(\nabla^2 \phi(\x_t)\Big)}{4k} 
\quad\mbox{and} \quad
\norm{S_t}_2 
\le
\left(\frac{5}{4} + \frac{d}{4k}\right)^{1/2},
\end{align*}
with a probability at least $1-\delta/T$.

For notation convenience, we denote $C_1 = 5\norm{ \nabla\phi(\x_t) }_2 + \frac{\tr(\nabla^2 \phi(\x_t))}{k}$ and $C_2 = \sqrt{\frac{5}{4} + \frac{d}{4k}}$.
We have
\begin{align*}
\phi(\x_{t+1}) - \phi(\x^*)
\stackrel{\eqref{eq:sas1}}{\leq}&
\phi(\x_t) - \phi(\x^*) - \eta_t \left(\frac{1}{2} - \frac{C_1}{4} \eta_t\right) \norm{S_t^\top \nabla \phi(\x_t)}^2 \\
+& \frac{4\mL\eta_t^3C_2^3}{3}\norm{S_t^\top \nabla\phi(\x_t)}^2\cdot \norm{\nabla \phi(\x_t)}\\
+&\eta_t\left(1 + 8C_1\eta_t\right)\left(\norm{\bw_t}^2 + \norm{\bv_t}^2\right)
+ \frac{8C_2^3\mL\eta_t^3}{3}\left(\norm{\bw_t}^3 + \norm{\bv_t}^3\right).
\end{align*}
Furthermore, if we have
\begin{align*}
\Big(\phi(\x_t) - \phi(\x^*)\Big)^{1/2} \le \frac{3\sqrt{2}C_1}{32\sqrt{L}\mL\eta_t C_2^3},
\end{align*}
then, we can obtain 
\begin{align*}
	\frac{4\mL\eta_t C_2^3}{3} \norm{\nabla \phi(\x_t)} 
	\le
	\frac{4\sqrt{2L}\mL\eta_t C_2^3}{3} \Big(\phi(\x_t) - \phi(\x^*)\Big)^{1/2} \le \frac{C_1}{4}.
\end{align*}

Combining above results, we have
\begin{align*}
\phi(\x_{t+1}) - \phi(\x^*)
\le&
\phi(\x_t) - \phi(\x^*) - \eta_t \left(\frac{1}{2} - \frac{C_1}{2} \eta_t\right) \norm{S_t^\top \nabla \phi(\x_t)}^2 \\
+&\eta_t\left(1 + 8C_1\eta_t\right)\left(\norm{\bw_t}^2 + \norm{\bv_t}^2\right)
+ \frac{8C_2^3\mL\eta_t^3}{3}\left(\norm{\bw_t}^3 + \norm{\bv_t}^3\right)\\
=&
\phi(\x_t) - \phi(\x^*) - \frac{1}{8C_1} \norm{S_t^\top \nabla \phi(\x_t)}^2  
+\frac{5}{2C_1}\left(\norm{\bw_t}^2 + \norm{\bv_t}^2\right)
+ \frac{C_2^3\mL}{3C_1^3}\left(\norm{\bw_t}^3 + \norm{\bv_t}^3\right)\\
\stackrel{\eqref{eq:sas1}}{\leq}&
\phi(\x_t) - \phi(\x^*) - \frac{1}{16C_1}\norm{\nabla \phi(\x_t)}^2 
+\frac{5}{2C_1}\left(\norm{\bw_t}^2 + \norm{\bv_t}^2\right)
+ \frac{C_2^3\mL}{3C_1^3}\left(\norm{\bw_t}^3 + \norm{\bv_t}^3\right)\\
\le&
\left(1 - \frac{\mu}{8C_1}\right)\Big(\phi(\x_t) - \phi(\x^*)\Big) 
+ \frac{5}{2C_1}\left(\norm{\bw_t}^2 + \norm{\bv_t}^2\right)
+ \frac{C_2^3\mL}{3C_1^3}\left(\norm{\bw_t}^3 + \norm{\bv_t}^3\right),
\end{align*}
where the last inequality is because of Lemma~\ref{lem:mu}.

By Eq.~\eqref{eq:g_prop1} and Lemma \ref{lem:bound_S_F_norm}, with probability at least $1-\delta/T$, we have
\begin{align*}
\norm{\bw_t}^2 + \norm{\bv_t}^2 
\le\frac{C_3^3\mL^2\alpha^4}{36} + \frac{\ell\sigma^2}{\alpha^2} \quad\mbox{ and }\quad \norm{\bw_t}^3 + \norm{\bv_t}^3 
\le\frac{C_3^{9/2}\mL^3\alpha^6}{6^3} + \frac{\ell^{3/2}\sigma^3}{\alpha^3}.
\end{align*}
Thus,
\begin{align*}
\phi(\x_{t+1}) - \phi(\x^*)
\leq
\left(1 - \frac{\mu}{8C_1}\right)\Big(\phi(\x_t) - \phi(\x^*)\Big) 
+ \frac{5}{2C_1} \cdot \left( \frac{C_3^3\mL^2\alpha^4}{36} + \frac{\ell\sigma^2}{\alpha^2} \right)
+ \frac{C_2^3\mL}{3C_1^3}\left( \frac{C_3^{9/2}\mL^3\alpha^6}{6^3} + \frac{\ell^{3/2}\sigma^3}{\alpha^3} \right),
\end{align*}
with probability at least $1-\delta/T$.
Rearranging above equation by Lemma~\ref{lem:dd}, we can obtain
\begin{align*}
\phi(\x_{t+1}) - \phi(\x^*) - C_4 \le \left(1 - \frac{\mu}{8C_1}\right)\cdot \Big(\phi(\x_t) - \phi(\x^*) - C_4\Big),
\end{align*}
with probability at least $1-\delta/T$.

\section{Proof of Theorems and Corollaries}

\subsection{Proof of Theorem \ref{thm:main}}
We have
\begin{align*}
\phi(\x_{t+1}) - \phi(\x^*)\stackrel{\eqref{eq:phi_dec}}{\leq}&
\phi(\x_t) - \phi(\x^*) - \frac{\eta}{2} \norm{S_t^\top \nabla\phi(\x_t)}^2 + \eta^2 \norm{S_t^\top A S_t}_2 \cdot \norm{S_t^\top \nabla \phi(\x_t)}^2 \\
+& \frac{\eta}{2} \norm{\bv_t}^2 + \eta^2 \norm{S_t^\top A S_t}_2 \cdot \norm{\bv_t}^2\\
\stackrel{\eqref{eq:sas1}}{\le}&
\phi(\x_t) - \phi(\x^*) - \frac{\eta}{2} \norm{S_t^\top \nabla\phi(\x_t)}^2 
+ \eta^2 \left(\frac{5\norm{A}_2}{4} + \frac{\tr(A)}{4k} \right) \cdot \norm{S_t^\top \nabla \phi(\x_t)}^2\\
+& \frac{\eta}{2} \norm{\bv_t}^2 + \eta^2 \left(\frac{5\norm{A}_2}{4} + \frac{\tr(A)}{4k} \right) \cdot \norm{\bv_t}^2\ \text{(with probability $1-\frac{\delta}{T}$)}\\
=&
\phi(\x_t) - \phi(\x^*) - \frac{1}{4}\cdot\left(5 \norm{A}_2 + \frac{\tr(A)}{k}\right)^{-1} \norm{S_t^\top \nabla\phi(\x_t)}^2 
+ \frac{3}{4}\left(5 \norm{A}_2 + \frac{\tr(A)}{k}\right)^{-1} \norm{\bv_t}^2\\
\stackrel{\eqref{eq:sas1}\eqref{eq:g_prop}}{\le}&
\phi(\x_t) - \phi(\x^*) - \frac{1}{8}\left(5 \norm{A}_2 + \frac{\tr(A)}{k}\right)^{-1} \norm{\nabla \phi(\x_t)}^2 + \frac{3\ell\sigma^2}{4\alpha^2} \left(5 \norm{A}_2 + \frac{\tr(A)}{k}\right)^{-1},
\end{align*}
where the first equality is because of $\eta = \left(5\norm{A}_2 + \frac{\tr(A)}{k}\right)^{-1}$.

By Assumption~\ref{ass:mu} and replacing $x, y$ with $\x^*$ and $\x_t$ respectively, then we have
\begin{align*}
	\phi(\x_t) - \phi(\x^*) \le \frac{1}{2\mu} \norm{\nabla \phi(\x_t)}^2.
\end{align*}

Therefore, with probability $1-\frac{\delta}{T}$
\begin{align*}
\phi(\x_{t+1}) - \phi(\x^*)
\le 
\left(1 - \frac{\mu}{4} \cdot \left(5 \norm{A}_2 + \frac{\tr(A)}{k}\right)^{-1}\right) \cdot \Big(\phi(\x_t) - \phi(\x^*)\Big) + \frac{3\ell\sigma^2}{4\alpha^2} \left(5 \norm{A}_2 + \frac{\tr(A)}{k}\right)^{-1}.
\end{align*}
By rearranging the above equations, we can obtain 
\begin{align*}
\phi(\x_{t+1}) - \phi(\x^*) - \frac{3\ell\sigma^2}{\mu \alpha^2}
\leq
\left(1 - \frac{\mu}{4} \cdot \left(5 \norm{A}_2 + \frac{\tr(A)}{k}\right)^{-1}\right) \cdot \left(\phi(\x_t) - \phi(\x^*) - \frac{3\ell\sigma^2}{\mu \alpha^2}\right),
\end{align*}
with probability $1-\frac{\delta}{T}$.

\subsection{Proof of Corollary \ref{cor:main}}
First, by the presentations in Section~\ref{subsec:ske}, choosing the sketching matrix $S_t\in\RR^{d\times \ell}$ with $\ell = \cO\left(k + \log\frac{T}{\delta}\right)$ is a $(\frac{1}{4}, k, \delta/T)$-oblivious sketching matrix for the matrix product.
Thus, by Theorem~\ref{thm:main}, we have,
\begin{align*}
\phi(\x_T) - \phi(\x^*) - \frac{3\ell\sigma^2}{\mu \alpha^2}
\leq&
\left(1 - \frac{\mu}{4} \cdot \left(5 \norm{A}_2 + \frac{\tr(A)}{k}\right)^{-1}\right) \cdot \left(\phi(\x_{T-1}) - \phi(\x^*) - \frac{3\ell\sigma^2}{\mu \alpha^2}\right)\\
\leq&
\left(1 - \frac{\mu}{4} \cdot \left(5 \norm{A}_2 + \frac{\tr(A)}{k}\right)^{-1}\right)^T \cdot \left(\phi(\x_0) - \phi(\x^*) - \frac{3\ell\sigma^2}{\mu \alpha^2}\right)\\
\leq&
\exp\left(- T \cdot\frac{\mu}{4} \cdot \left(5 \norm{A}_2 + \frac{\tr(A)}{k}\right)^{-1}\right)\left(\phi(\x_0) - \phi(\x^*) - \frac{3\ell\sigma^2}{\mu \alpha^2}\right),
\end{align*}
with probability at least  $1-{\delta}$.
Letting the right-hand side of the above equation equal to $\epsilon$, then it requires that
\begin{align*}
	T = \frac{20 \norm{A}_2 + \frac{4\tr(A)}{k}}{\mu}\log\frac{\left(\phi(\x_0) - \phi(\x^*) - \frac{3\ell\sigma^2}{\mu \alpha^2}\right)}{\epsilon} 
	= \cO\left(\left(\frac{\norm{A}_2}{\mu} + \frac{\tr(A)}{\mu\ell }\right)\log\frac{1}{\epsilon}\right),
\end{align*}
where the last equality is because the assumption $k > \log\frac{T}{\delta}$.

Accordingly, the total query complexity is
\begin{align*}
	Q = T\times\ell = \cO\left(\left(\frac{\ell\norm{A}_2}{\mu} + \frac{\tr(A)}{\mu}\right)\log\frac{1}{\epsilon}\log\left(\left(\frac{\norm{A}_2}{\mu} + \frac{\tr(A)}{\mu\ell }\right)\log\frac{1}{\epsilon}\right)\right)=\Tilde{\cO}\left(\left(\frac{\ell\norm{A}_2}{\mu} + \frac{\tr(A)}{\mu}\right)\log\frac{1}{\epsilon}\right).
\end{align*}

\subsection{Proof of Theorem \ref{thm:main1}}

	By the update rule and Assumption~\ref{ass:phi}, we have
	\begin{align*}
		\phi(\x_{t+1}) 
		=&
		\phi(\x_t) - \eta_t \dotprod{\nabla \phi(\x_t), \tg(\x_t)} + \frac{\eta_t^2}{2} \norm{\tg(\x_t)}_A^2\\
		\stackrel{\eqref{eq:tg_prop}}{=}&
		\phi(\x_t) - \eta_t\dotprod{\nabla \phi(\x_t), H_t^{-1/2}S_tS_t^\top H_t^{-1/2} \nabla \phi(\x_t) + H_t^{-1/2}S_t\tv_t} \\
		 + &\frac{\eta_t^2}{2} \norm{ H_t^{-1/2}S_tS_t^\top H_t^{-1/2} \phi(\x_t) + H_t^{-1/2}S_t\tv_t }^2_A\\
		\le&
		\phi(\x_t) - \eta_t \norm{S_t^\top H_t^{-1/2} \nabla \phi(\x_t)}^2  - \eta_t \dotprod{\nabla \phi(\x_t), H_t^{-1/2}S_t\tv_t}\\
		 +& \eta^2 \norm{H_t^{-1/2}S_tS_t^\top H_t^{-1/2} \phi(\x_t)}_A^2 + \eta^2 \norm{ H_t^{-1/2}S_t\tv_t }_A^2\\
		\le&
		\phi(\x_t) - \frac{\eta_t}{2} \norm{S_t^\top H_t^{-1/2} \nabla \phi(\x_t)}^2 + \frac{\eta_t}{2} \norm{\tv_t}^2 +  \eta^2 \norm{H_t^{-1/2}S_tS_t^\top H_t^{-1/2} \phi(\x_t)}_A^2 + \eta^2 \norm{ H_t^{-1/2}S_t\tv_t }_A^2\\
		\le&
		\phi(\x_t) - \frac{\eta_t}{2} \norm{S_t^\top H_t^{-1/2} \nabla \phi(\x_t)}^2 + \eta_t^2 \norm{S_t^\top H_t^{-1/2} A H_t^{-1/2} S_t}_2 \cdot \norm{S_t^\top H_t^{-1/2} \nabla \phi(\x_t)}^2 \\
	+ &\eta_t^2\norm{S_t^\top H_t^{-1/2} A H_t^{-1/2} S_t}_2 \cdot \norm{\tv_t}^2 + \frac{\eta_t}{2} \norm{\tv_t}^2.
	\end{align*}
	where the second inequality is because of Cauchy's inequality and $2ab \le a^2 + b^2$.
	
	Replacing $A$ and $\nabla \phi(\x_t)$ with $H_t^{-1/2} A H^{-1/2}$ and $H_t^{-1/2}\nabla \phi(\x_t)$ respectively in Lemma~\ref{lem:ske}, we can obtain that
	\begin{align*}
		\norm{S_t^\top H_t^{-1/2} A H_t^{-1/2} S_t}_2 
		\leq
		\frac{5\norm{ H_t^{-1/2} A H_t^{-1/2} }_2}{4} + \frac{\tr(H_t^{-1} A)}{4k},
	\end{align*} 
	and
	\begin{align*}
		\norm{S_t^\top H_t^{-1/2}\nabla \phi(\x_t)}^2 \ge \frac{1}{2} \norm{ H_t^{-1/2}\nabla \phi(\x_t) }^2
	\end{align*}
    with probability $1-\frac{\delta}{T}$.
	Combining above results with the step size $\eta_t = \left(  5\norm{ H_t^{-1/2} A H_t^{-1/2} }_2 + \frac{\tr(H_t^{-1} A)}{k}\right)^{-1}$, we can obtain that
	\begin{align*}
	&\phi(\x_{t+1}) -  \phi(\x^*)\\
	\leq&
	\phi(\x_t) - \phi(\x^*) - \eta_t \left(\frac{1}{2} - \eta_t \left(  \frac{5\norm{ H_t^{-1/2} A H_t^{-1/2} }_2}{4} + \frac{\tr(H_t^{-1} A)}{4k}\right)\right) \cdot \norm{S_t^\top H_t^{-1/2} \nabla \phi(\x_t)}^2 \\
	+&  \eta_t\left(\frac{1}{2} + \eta_t \left(  \frac{5\norm{ H_t^{-1/2} A H_t^{-1/2} }_2}{4} + \frac{\tr(H_t^{-1} A)}{4k}\right)\right)  \cdot \norm{\tv_t}^2\\
	=&
	\phi(\x_t) - \phi(\x^*) - \frac{1}{4}\left( \norm{ H_t^{-1/2} A H_t^{-1/2} }_2 + \frac{\tr(H_t^{-1} A)}{k} \right) \cdot \norm{S_t^\top H_t^{-1/2} \nabla \phi(\x_t)}^2 \\
	+& \frac{3}{4} \left( \norm{ H_t^{-1/2} A H_t^{-1/2} }_2 + \frac{\tr(H_t^{-1} A)}{k} \right) \cdot \norm{\tv_t}^2\\
	\le&
	\phi(\x_t) - \phi(\x^*) - \frac{1}{8}\left( \norm{ H_t^{-1/2} A H_t^{-1/2} }_2 + \frac{\tr(H_t^{-1} A)}{k} \right)^{-1} \cdot \norm{ H_t^{-1/2} \nabla \phi(\x_t)}^2 \\
	+& \frac{3}{4} \left( \norm{ H_t^{-1/2} A H_t^{-1/2} }_2 + \frac{\tr(H_t^{-1} A)}{k} \right)^{-1} \cdot \norm{\tv_t}^2\\
	\stackrel{\eqref{eq:rho}}{\leq}&
	\left(1 - \frac{\rho_t}{4} \left( \norm{ H_t^{-1/2} A H_t^{-1/2} }_2 + \frac{\tr(H_t^{-1} A)}{k} \right)^{-1} \right)\cdot \Big(\phi(\x_t) - \phi(\x^*)\Big)  \\
	+&  \left( \norm{ H_t^{-1/2} A H_t^{-1/2} }_2 + \frac{\tr(H_t^{-1} A)}{k} \right)^{-1} \cdot \frac{3\ell\sigma^2}{4\alpha^2},\ \ \text{with probability $1-\frac{\delta}{T}$,}
	\end{align*}
	where the first equality is because of the condition on the step size $\eta_t$.
	
	Rearrange the above equation, we can obtain that
	\begin{align*}
	\phi(\x_{t+1}) -  \phi(\x^*) - \frac{3\ell\sigma^2}{\rho_t \alpha^2} 
	\leq
	\left(1 - \frac{\rho_t}{4} \left( \norm{ H_t^{-1/2} A H_t^{-1/2} }_2 + \frac{\tr(H_t^{-1} A)}{k} \right)^{-1} \right)\cdot \left(\phi(\x_t) - \phi(\x^*) - \frac{3\ell\sigma^2}{\rho_t \alpha^2}\right)
	\end{align*}
    with probability $1-\frac{\delta}{T}$.

\subsection{Proof of Corollary \ref{cor:main1}}

	$S_t\in\RR^{d\times \ell}$ is the $(\frac{1}{4},k,\frac{\delta}{T})$-oblivious sketching matrix with $\ell = \cO\left(k + \log\frac{T}{\delta}\right)$ for the matrix product.
	Thus, by Theorem~\ref{thm:main1}, we have,
	\begin{align*}
		\phi(\x_T) - \phi(\x^*)
		\leq
		\exp\left(- T \cdot\min_t\left\{\frac{\rho_t}{4} \cdot \left(5 \norm{H_t^{-1/2}AH_t^{-1/2}}_2 + \frac{\tr(H_t^{-1}A)}{k}\right)^{-1}\right\}\right)\left(\phi(\x_0) - \phi(\x^*) - \frac{3\ell\sigma^2}{\mu \alpha^2}\right) +  \frac{3\ell\sigma^2}{\rho \alpha^2},
	\end{align*}
    with probability $1-\delta$.
	Letting the first term of the right-hand side of the above equation equal to $\epsilon$ and using the fact $\ell = \cO(k)$, then it requires that
	\begin{align*}
		T 
		=\cO\left(\max_t\left\{\rho_t^{-1}\left( \norm{ H_t^{-1/2} A H_t^{-1/2} }_2 + \frac{\tr(H_t^{-1} A)}{\ell} \right)\right\}\log\frac{1}{\epsilon} \right).
	\end{align*}
	
	Accordingly, the total query complexity is
	\begin{align*}
		Q = T\times\ell = \Tilde{\cO}\left(\max_t\left\{\rho_t^{-1}\left( \norm{ H_t^{-1/2} A H_t^{-1/2} }_2 \cdot \ell + \tr(H_t^{-1} A) \right)\right\}\log\frac{1}{\epsilon} \right).
	\end{align*}

\subsection{Proof of Theorem \ref{thm:trace}}

First, using Taylor's expansion, we have
\begin{align*}
&\sum_{i=1}^{\ell} \frac{\phi(\x+\alpha \sti) + \phi(\x_t-\alpha \sti) - 2\phi(\x_t)}{ \alpha^2} \\
=&
\sum_{i=1}^{\ell} \frac{\phi(\x_t) + \alpha\dotprod{\nabla \phi(\x_t), \sti} + \frac{\alpha^2}{2} \norm{\sti}_{\nabla^2 \phi(\x_t)}^2 + \phi(\x_t) - \alpha\dotprod{\nabla \phi(\x_t), \sti} + \frac{\alpha^2}{2} \norm{\sti}_{\nabla^2 \phi(\x_t)}^2 - 2\phi(\x_t)}{\alpha^2} \\
+& \frac{D(\alpha) + D(-\alpha)}{\alpha^2}\\
=&
\sum_{i=1}^{\ell} (\sti)^\top \nabla^2 \phi(\x_t) \sti + \frac{D(\alpha, \sti) + D(-\alpha, \sti)}{\alpha^2},
\end{align*}
where we define that $D(\alpha, \sti) =  \phi(\x_t+\alpha \sti) - \left( \phi(\x_t) + \alpha\dotprod{\nabla \phi(\x_t), \sti} + \frac{\alpha^2}{2} \norm{\sti}_{\nabla^2 \phi(\x_t)}^2 \right)$.

Furthermore,
\begin{align*}
	\left|\frac{D(\alpha, \sti) + D(-\alpha, \sti)}{\alpha^2}\right| 
	\le \frac{|D(\alpha, \sti)|+|D(-\alpha, \sti)|}{\alpha^2} 
	\stackrel{\eqref{eq:mL}}{\leq} \frac{\mL\norm{\sti}^3}{3}\alpha. 
\end{align*}

Moreover, by Lemma~\ref{lem:trace}, we can obtain that 
\begin{equation*}
\left| \sum_{i=1}^{\ell} (\sti)^\top \nabla^2 \phi(\x_t) \sti - \tr\left(\nabla^2 \phi(\x_t)\right) \right| \le \varepsilon \cdot \tr\left(\nabla^2 \phi(\x_t)\right).
\end{equation*}
Combining with the definition of $f(\x_t)$ and Assumption~\ref{ass:zeta}, we have
\begin{align*}
	&\left|\tau(\x_t, S_t) - \tr\Big(\nabla^2 \phi(\x_t)\Big)\right| \\
	= &
	\left|\sum_{i=1}^{\ell} \frac{\phi(\x_t+\alpha \sti) + \phi(\x_t-\alpha \sti) - 2\phi(\x_t) + \zeta(\x_t+\alpha \sti)+ \zeta(\x_t-\alpha \sti)-2\zeta(\x_t)}{ \alpha^2} - \tr\Big(\nabla^2 \phi(\x_t)\Big) \right|\\
	\le&
	\left| \sum_{i=1}^{\ell} (\sti)^\top \nabla^2 \phi(\x_t) \sti -\tr\Big(\nabla^2 \phi(\x_t)\Big)   \right| + \sum_{i=1}^{\ell} \left| \frac{D(\alpha, \sti) + D(-\alpha, \sti)}{\alpha^2}    \right| \\
	+& \sum_{i=1}^{\ell} \left| \frac{\zeta(\x_t+\alpha \sti)+ \zeta(\x_t-\alpha \sti)-2\zeta(\x_t)}{ \alpha^2} \right| \\
	\le&
	\varepsilon\cdot \tr\Big(\nabla^2 \phi(\x_t)\Big) + \sum_{i=1}^{\ell}\frac{\mL\norm{\sti}^3}{3}\alpha + \frac{4\ell\sigma}{\alpha^2}.
\end{align*}

For the Rademacher and SRHT sketching matrix, it holds that $\norm{\sti} = \sqrt{d / \ell}$. 
Thus, we have $\sum_{i=1}^{\ell} \norm{\sti}^3 = d^{3/2}\ell^{-1/2}$.

For the Gaussian sketching matrix, we  first have
\begin{align*}
	\norm{\sti}^2 \le \frac{2d +3\log\frac{1}{\delta} }{\ell}.
\end{align*}
Thus, with a probability at least $1-\delta$,
\begin{align*}
\sum_{i=1}^\ell\norm{\sti}^{3/2} \le \left(2d + \log\frac{\ell}{\delta}\right)^{3/2} \ell^{-1/2}.
\end{align*}

Combining the above results, we can obtain the final result.



\section{Other Useful Lemmas}

\begin{lemma}[Theorem 2.1.10 of \cite{nesterov2003introductory}] \label{lem:mu}
	If $\phi(\x)$ is $\mu$-strongly convex and differentiable, then it holds that
	\begin{equation*}
		\phi(\y) \le \phi(\x) + \dotprod{\nabla \phi(\x), \y - \x} + \frac{1}{2\mu} \norm{\nabla \phi(\x) - \nabla \phi(\y)}^2.
	\end{equation*}
\end{lemma}


\begin{lemma}
Letting $A\in\RR^{d\times d}$ be positive definite, then it holds that
\begin{equation}\label{eq:tr_F}
\norm{A^{1/2}}_F^2 = \tr(A), \quad\mbox{and}\quad\norm{A^{1/2}}_2^2 = \norm{A}_2.
\end{equation}
\end{lemma}

\textbf{Proof}: Letting $A = U\Lambda U^\top$ be eigenvalue decomposition of $A$ with $\Lambda = \diag(\lambda_1,\dots,\lambda_d)$. 
	Then, 
\begin{align*}
	\norm{A^{1/2}}_F^2 =& \sum_{i=1}^{d} \left(\lambda_i^{1/2}\right)^2 = \sum_{i=1}^{d} \lambda_i = \tr(A),\\
	\norm{A^{1/2}}_2^2 =& \left(\lambda_{1}^{1/2}\right)^2 = \lambda_1 = \norm{A}_2.
\end{align*}

\begin{lemma}
	Let Assumption~\ref{ass:phi}-\ref{ass:mu} hold. The approximate Hessian $H_t$ satisfies Eq.~\eqref{eq:preceq}, then it holds that
	\begin{equation}\label{eq:rho}
	\norm{H_t^{-1/2} \nabla \phi(\x_t) }^2 \ge 2\rho_t \Big(\phi(\x_t) - \phi(\x^*)\Big).
	\end{equation}
\end{lemma}

\textbf{Proof}:	By the property of the convex quadratic function, we have
	\begin{align*}
		\x^* = \x - A^{-1} \nabla \phi(\x).
	\end{align*}
	Replacing $\y$ with above equation in Eq.~\eqref{eq:phi}, we have
	\begin{align*}
		\phi(\x^*) 
		= 
		\phi(\x) - \norm{\nabla \phi(\x)}_{A^{-1}}^2 + \frac{1}{2} \norm{A^{-1}\nabla \phi(\x)}_A^2
		= \phi(\x) - \frac{1}{2} \norm{\nabla \phi(x)}_{A^{-1}}^2.
	\end{align*}
	Thus, we have
	\begin{align*}
		\norm{H_t^{-1/2} \nabla \phi(\x_t) }^2 
		= 
		\norm{\nabla \phi(\x_t)}_{H_t^{-1}}^2 
		\stackrel{\eqref{eq:preceq}}{\geq}&
		\rho_t   \norm{\nabla \phi(\x_t)}_{A^{-1}}^2 
		= 2\rho_t \Big(\phi(\x_t) - \phi(\x^*)\Big).
	\end{align*}

\begin{lemma}\label{lem:bound_S_F_norm}
    If $S\in\RR^{d\times \ell}$ is an oblivious  $(\frac{1}{4}, k, \delta)$-random sketching matrix, then with a probability at least $1-\delta$, it holds that
    \begin{equation}
        \norm{S}_F^2 \leq \frac{5\ell}{4} + \frac{d \ell}{4k}.
    \end{equation}
\end{lemma}

\textbf{Proof}:
    By the property of an oblivious  $(\frac{1}{4}, k, \delta)$-random sketching matrix, we can obtain that 
    \begin{equation*}
        \norm{SS^\top - I_{d}}_2 \leq \frac{1}{4} \left(\norm{I}_2 + \frac{\norm{I_{d}}_F^2}{k}\right) = \frac{1}{4}\left(1 + \frac{d}{k}\right),
    \end{equation*}
    which implies that
    \begin{equation*}
        \norm{S}_2^2 \leq \frac{5}{4} + \frac{d}{4k}.
    \end{equation*}
    Combining with the fact that $\norm{S}_F^2 \le \ell \norm{S}_2^2$, we can obtain that
    \begin{equation*}
        \norm{S}_F^2 \leq \frac{5\ell}{4} + \frac{d \ell}{4k}.
    \end{equation*}

\begin{remark}
    Above lemma holds for all oblivious  $(\frac{1}{4}, k, \delta)$-random sketching matrix. 
    However, for the Radamecher sketching matrix and SRHT, one can easily check that $\norm{S}_F^2 = d$. Therefore, we denote that $C_3=d$ for Radamecher sketching matrix and SRHT and $C_3= \frac{5\ell}{4} + \frac{d \ell}{4k}$ for other sketching matrices.
\end{remark}

\begin{lemma}
	\label{lem:dd}
	Letting non-negative sequence $\{\Delta_t\}$ satisfy 
$\Delta_{t+1} \leq (1 - \beta) \Delta_t + c$ with $0\leq \beta \leq 1$ and $c\geq 0$, then it holds that $\Delta_{t+1} - \frac{c}{\beta} \leq (1-\beta)\left(\Delta_t - \frac{c}{\beta}\right)$.
\end{lemma}

\textbf{Proof}:
	We assume that $\{\Delta_t\}$ satisfies that $ \Delta_{t+1} - c' \leq \left(1 - \beta\right)\left(\Delta_t - c'\right) $ with $c'\geq 0$. 
	Then, it implies that 
	\begin{equation*}
		\Delta_{t+1} \leq  (1 - \beta) \Delta_t + c' - (1-\beta)c' = (1 - \beta) \Delta_t + \beta c'.
	\end{equation*}
	Thus, we can obtain that $c' = \beta^{-1} c$ which concludes the proof.

\begin{lemma}[Corollary 2 of \cite{cortinovis2021randomized}]\label{lem:trace}
	Let $A \in \RR^{d\times d}$ be positive definite, $\delta \in (0, 1/2]$, $\varepsilon\in(0,1)$, $\ell \in \mathbb{N}$. Let $\tau_\ell(A) = \frac{1}{\ell}\sum_{i=1}^{\ell} \s_i^\top A \s_i$. To achieve the following property  holding with a probability of at least $1-\delta$,
	\begin{equation*}
		\Big|\tau_\ell(A) - \tr(A)\Big| \le \varepsilon \cdot \tr(A),
	\end{equation*}
	then, 
	\begin{enumerate}
		\item if $\s_i$ is a random standard Gaussian vector with zero mean, $\ell = 8\varepsilon^{-2}\frac{\norm{A}_2}{\tr(A)}\log\frac{2}{\delta}$ is sufficient.
		\item  if $\s_i \in\RR^d$ is a random Rademacher vector, $\ell = 8\varepsilon^{-2}(1+\varepsilon)\frac{\norm{A}_2}{\tr(A)}\log\frac{2}{\delta}$ is sufficient.
		\item if $\s_i \in\RR^d$ is a column of SRHT, $\ell = \left(1 + \sqrt{8\ln\frac{2n}{\delta}}\right)^2\varepsilon^{-2}\log\frac{2}{\delta}$ is sufficient.
	\end{enumerate}
\end{lemma}

\begin{lemma}[Lemma 2 of \cite{meyer2021hutch}]\label{lem:trace1}
	Let $A \in \RR^{d\times d}$ be positive definite, $\delta \in (0, 1/2]$, $\ell \in \mathbb{N}$. Let $\tau_\ell(A) = \frac{1}{\ell}\sum_{i=1}^{\ell} \s_i^\top A \s_i$, where $\s_i \in\RR^d$ is mean 0, i.i.d. sub-Gaussian random variable with constant sub-Gaussian parameter $C$.
	For fixed constants $c, C$ which $c$ only depends on $C$, if $\ell \ge \cO\left(c \log\frac{1}{\delta}\right)$, then with a probability $1-\delta$, 
	\begin{equation}\label{eq:trace}
		|\tau_\ell(A) - \tr(A)| \le C\sqrt{\frac{\log\frac{1}{\delta}}{\ell}} \norm{A}_F.
	\end{equation}
\end{lemma}

\end{APPENDICES}



\end{document}